\documentclass[12pt]{amsart}
\usepackage{bm}

 \newtheorem{thm}{Theorem}[section]
 \newtheorem{cor}[thm]{Corollary}
 \newtheorem{lem}[thm]{Lemma}
 \newtheorem{prop}[thm]{Proposition}
 \theoremstyle{definition}
 \newtheorem{defn}[thm]{Definition}
 \theoremstyle{remark}
 \newtheorem{rem}[thm]{Remark}
 \numberwithin{equation}{section}

 \newcommand{\Real}{\mathbb{R}}

 \newcommand{\dds}{\frac{d}{ds}}
 \newcommand{\ddt}{\frac{d}{dt}}

 \newcommand{\ent}[2]{\textbf{Ent}(#1|#2)}
 \newcommand{\ric}{\mathfrak{Ric}}
\begin{document}

\title[Ricci Curvature for Contact 3-Manifolds]
{Generalized Ricci Curvature Bounds for Three Dimensional Contact Subriemannian Manifolds}

\author{Andrei Agrachev}
\email{agrachev@sissa.it}
\address{International School for Advanced Studies
via Beirut 4 - 34014, Trieste, Italy and Steklov Mathematical Institute, ul. Gubkina 8,
Moscow, 119991 Russia}

\author{Paul W.Y. Lee}
\email{wylee@math.cuhk.edu.hk}
\address{Room 216, Lady Shaw Building, The Chinese University of Hong Kong, Shatin, Hong Kong}

\date{\today}

\thanks{The first author was partially supported by PRIN and the second author was supported by NSERC postgraduate scholarship and postdoctoral fellowship.}

\begin{abstract}
Measure contraction property is one of the possible generalizations of Ricci curvature bound to more general metric measure spaces. In this paper, we discover necessary and sufficient conditions for a three dimensional contact subriemannian manifold to satisfy this property.
\end{abstract}

\maketitle

\section{Introduction}

In the past few years, several connections between the optimal transportation problems and curvature of Riemannian manifolds were found. One of them is the use of optimal transportation for an alternative definition of Ricci curvature lower bound developed in a series of papers \cite{OtVi,CoMcSc,StVo}. Based on the ideas in these papers, a generalization of Ricci curvature lower bound for general metric measure spaces, called curvature-dimension condition, is introduced in \cite{LoVi1,LoVi2,St1,St2}  (see section \ref{optimalRicci} for a quick overview of these results). Recently the case of a Finsler manifold was studied in \cite{Oh2} and the results are very similar to that of the Riemannian case due to strict convexity of the corresponding Hamiltonian.

The situation changes dramatically in the case of subriemannian manifolds. The reason is that the class of metric spaces we are dealing with have Hausdorff dimensions strictly greater than their topological dimensions. Therefore, the interplay between the metrics and the measures of these spaces should be significantly different from that of the Riemannian or Finsler case. One particular case of subriemannian manifolds, the Heisenberg group, is studied in \cite{Ju}. In this case the space does not satisfy any curvature-dimension condition mentioned above (however, see \cite{BaGa1,BaBoGa,BaGa2} for a different definition of curvature-dimension condition in the subriemannian setting). Instead it satisfies a weaker condition, called measure contraction property, introduced in \cite{St2,Oh1} (see Section \ref{optimalRicci} for the definition). Like the curvature-dimension condition, measure contraction property is a generalization of Ricci curvature lower bound on Riemannian manifolds. However, it is a weaker condition for general metric measure spaces.

The approach used by \cite{Ju} relies on the complete integrability of the subriemannian geodesic flow on the Heisenberg group. Because of this, the changes in the measure along the geodesic flow can be written down explicitly in this case, which is not possible for subriemannian manifolds in general.

The goal of this paper is to study a subriemannian version of the measure contraction property for three dimensional contact subriemannian manifolds under certain curvature conditions. This study uses a subriemannian generalization of the classical Riemannian curvature. The generalized Ricci curvature was introduced by the first author in the 90s for some special cases (including the three dimensional contact subriemannian structures), and in full generality by the first author and I.~Zelenko (see \cite{AgZe}). Later C.-B. Li and I. Zelenko found a complete system of curvature invariants (see \cite{LiZe1,LiZe2}). To state some interesting consequences of the main result in this paper, let us first give a brief introduction to the curvature invariants (see Section \ref{gencur} for a more detailed discussion of these invariants).

Let $e^{t\vec H}$ be the subriemannian geodesic flow defined on the cotangent bundle $T^*M$ of a manifold $M$ and let $\alpha$ be in $T^*M$. In a similar spirit of the Fr\'enet-–Serret frame, one can find a special moving frame along the trajectory $t\mapsto e^{t\vec H}(\alpha)$. The main property of this frame is that it satisfies certain first order equations when pulled back to the tangent space $T_\alpha T^*M$ at the point $\alpha$ by the geodesic flow $e^{t\vec H}$. The pulled back frame is called the canonical Darboux frame.

In the Riemannian case, the canonical Darboux frame
\[
\{e_1(t),...,e_n(t),f_1(t),...,f_n(t)\}
\]
satisfies the following equations which is the Jacobi field equation (up to certain identifications of tangent and cotangent spaces)
\[
\dot e_i(t)=f_i(t),\quad \dot f_i(t)=-R^{ij}_\alpha(t)e_j(t).
\]
The matrix $R_\alpha:=R_\alpha(0)$ with $ij$-th entries given by $R^{ij}_\alpha(0)$ above is the Riemannian curvature operator (again up to certain identifications).

In the three dimensional contact subriemannian case, the canonical Darboux frame
\[
\{e_1(t),e_2(t),e_3(t),f_1(t),f_2(t),f_3(t)\}
\]
satisfies the following equations instead
\[
\begin{array}{ll}
    \dot e_1(t)=f_1(t), \quad \dot f_1(t)=-R^{11}_\alpha(t)e_1(t)-f_2(t),\\
    \dot e_2(t)=e_1(t), \quad     \dot f_2(t)=-R^{22}_\alpha(t)e_2(t),\\
    \dot e_3(t)=f_3(t), \quad    \dot f_3(t)=0.
\end{array}
\]
Therefore, $R_\alpha=\left(\begin{array}{ccc}
R^{11}_\alpha(0) & 0 & 0\\
0 & R^{22}_\alpha(0) & 0\\
0 & 0 & 0
\end{array}\right)$ is a natural generalization of the Riemannian curvature.

In this paper, we introduce a new generalized measure contraction property $\mathcal{MCP}(K;2,3)$ (see Section \ref{GMCP} for the definition and its motivation). One of the main results (Theorem \ref{mainSasa}) gives necessary and sufficient conditions on the curvature $R$ for a class of three dimensional contact subriemannian manifolds, called Sasakian manifolds, to satisfy this new measure contraction property. Our generalized measure contraction property $\mathcal{MCP}(0;2,3)$ coincides with the old condition $MCP(0,5)$ (see Section \ref{optimalRicci} for the definition of $MCP(0,5)$). As a result of this and Theorem \ref{mainSasa}, the following theorem holds. In particular, it generalizes the result in \cite{Ju} for the Heisenberg group.

\begin{thm}\label{SasaMCP}(Measure Contraction Property)
Assume that the three dimensional contact subriemannian manifold $M$ is Sasakian. If $R^{11}_\alpha\geq 0$ for all $\alpha$ in the cotangent bundle $T^*M$, then the metric measure space $(M,d,\eta)$ satisfies the measure contraction property $MCP(0,5)$, where $d$ is the natural subriemannian structure defined on $M$ and $\eta$ is the corresponding Popp measure (see Section 3 for the precise definitions).
\end{thm}

Several interesting consequences also follow from Theorem \ref{SasaMCP} (see Section \ref{main} for the detail). They include:

\begin{itemize}
\item Volume doubling property
\item Local Poincar\'e inequality
\item Harnack inequality for harmonic functions of sub-Laplacian
\item Liouville property of sub-Laplacian
\end{itemize}

The method used in the proof of Theorem \ref{mainSasa} also apply to three dimensional contact subriemannian manifolds which are not necessarily Sasakian. In the second main result (Theorem \ref{morethm}), we apply it to any three dimensional compact contact subriemannian manifolds and give estimates of the measure contractions for these subriemannian manifolds.

The structure of this paper is as follows. In Section \ref{subgeo}, we give several basic notions on subriemannian  geometry necessary for the present work. In Section \ref{contact}, we give the definition and the properties of contact subriemannian manifolds. A special class of examples of contact subriemannian manifolds, called Sasakian manifolds, is  introduced here as well. Sasakian manifolds serve as examples to the main result of this paper. In Section \ref{OptimalIntro}, we recall the definition and some basic results on the optimal transportation problem. In Section \ref{optimalRicci},  we give a brief overview on how the optimal transportation problem gives rise to the curvature-dimension condition and the measure contraction property. We also give the motivation of the present work in this section. In Section \ref{gencur}, we recall and specialize the recent result of \cite{LiZe1,LiZe2} on the curvature type invariants of subriemannian manifolds to the three dimensional contact case. We also give explicit formulas for these invariants. In section \ref{GMCP}, we give the definition of the new generalized measure contraction property. We will also motivate its definition by considering how measures contract in the Sasakian version of space forms. In Section \ref{main}, we state the main theorem (Theorem \ref{mainSasa}) and its consequences. The main theorem gives necessary and sufficient conditions on when a three dimensional Sasakian manifold equipped with the natural subriemannian structure and Popp's measure satisfies the generalized measure contraction property $\mathcal{MCP}(K;2,3)$. (see Definition \ref{genMCP} below). In particular, $\mathcal{MCP}(0;2,3)$ coincides with the old measure contraction property $MCP(0,5)$. As a consequence, these spaces satisfy the volume doubling property, the local Poincar\'e inequality, the Harnack inequality and the Liouville property for harmonic functions of the sub-Laplacian. In Section \ref{more}, we give the measure contraction estimates for three dimensional compact contact subriemannian manifolds. The proofs of all the results of this paper are given in the rest of the sections.

\smallskip
\smallskip

Soon after we posted the first version of this paper on arXiv, there are some very interesting related works emerge. They are also related to Ricci curvature type condition and its consequences to subriemannian geometry and PDEs (\cite{BaGa1,LiZe3,BaBoGa,BaGa2,AgLe2}) . Among them, \cite{BaGa1,BaBoGa,BaGa2} uses an approach very different from ours. It would be very interesting to establish connections between the two approaches.

We would also like to mention a few closely related works appeared earlier which were pointed out by the referees. Comparison type results for contact manifolds using a different approach were considered in \cite{Ru,Hu}. In particular, a Bonnet-Myer theorem was proved there. In the Sasakian case, a volume comparison theorem were also considered in \cite{ChYa} using a different approach and some of the results in Theorem \ref{mainSasa} can be done using the Wronskian comparison theorem proved there. However, unlike the Riemannian case, the volume comparison theorem and measure contraction properties are different in the subriemannian case (see \cite{AgLe2}). Finally, Jacobi fields on Sasakian manifolds in higher dimensions were also considered earlier in \cite{BaDr}.

\smallskip

\section*{Acknowledgment}
We thank Igor Zelenko and Cheng-Bo Li for very interesting and stimulating discussions. This work is part of the PhD thesis of the second author. He  would like to express deep gratitude to his supervisor, Boris Khesin, for his continuous support. He would like to thank Professor Karl-Theodor Sturm for the fruitful discussions. He is also grateful to SISSA for their kind hospitality where part of this work is done. Finally, we would also like to thank the referees for providing many constructive comments.

\bigskip
\bigskip

\begin{tabular}{|l|r|}
  \hline
  \multicolumn{2}{|c|}{A table of notations} \\
  \hline
  $M$ & a metric space or a manifold  \\ \hline
    $\left<\cdot,\cdot\right>$, $\left|\cdot\right|$ & a subriemannian metric and its norm  \\ \hline
  $\Delta$ & a distribution on $M$ \\ \hline
  $\eta$ & the Popp's measure \\ \hline
  $\mu$, $\mu_0,\mu_1$, $\mu_t$ & measures on $M$ \\ \hline
  $\Pi$ & a measure on $M\times M$ \\ \hline
  $d$ & a distance function on $M$ \\ \hline
  $U$ & a Borel set in $M$  \\ \hline
  $H$ & subriemannian Hamiltonian \\ \hline
  $e^{t\vec H}$ & subriemannian geodesic flow\\ \hline
  $e_i(t),f_i(t)$ & canonical Darboux frame \\ \hline
  $R^{ij}_\alpha(t)$ & curvature invariants \\ \hline
  $v$ & tangent vectors \\ \hline
  $\alpha$ & covectors \\ \hline
  $\alpha_0$ & contact form \\ \hline
  $v_0$ & the Reeb field \\ \hline
  $v_1,v_2$ & subriemannian orthonormal basis \\ \hline
  $\alpha_0,\alpha_1,\alpha_2$ & dual basis of $v_0,v_1,v_2$\\ \hline
  $\theta$ & tautological 1-form on $T^*M$ \\ \hline
  $\omega$ & standard symplectic 2-form on $T^*M$ \\ \hline
  $X$ & a tangent vector in $TT^*M$\\\hline
  $\mathcal L$ & Lie derivative \\\hline
  $\nabla_H$ & horizontal gradient \\\hline
  $\Delta_H$ & sub-Laplacian \\\hline
\end{tabular}

\smallskip

\section{Subriemannian Manifolds and Their Geodesics}\label{subgeo}

In this section, we recall several basic notions in subriemannian geometry. For a detail discussion of various topics, see \cite{Mo}.

Recall that a Riemannian manifold is a manifold $M$ together with a fibrewise inner product defined on the tangent bundle $TM$. The length of a curve is defined by this inner product and the Riemannian distance between two points is the length of the shortest curve connecting them. For a subriemannian manifold the fibrewise inner product is defined on a family of subspaces $\Delta$ inside the tangent bundle $TM$. Therefore, the notion of length can only be defined for curves which are tangent to this family $\Delta$. These curves are called horizontal curves and the subriemannian distance between two points is the length of the shortest horizontal curve connecting them.

More precisely, a subriemannian manifold is a triple $(M,\Delta,\left<\cdot,\cdot\right>)$, where $M$ is a smooth manifold, $\Delta$ is a distribution (a vector subbundle $\Delta$ of the tangent bundle $TM$ of the manifold $M$), and $\left<\cdot,\cdot\right>$ is a fibrewise inner product defined on the distribution $\Delta$. The inner product $\left<\cdot,\cdot\right>$ is also called a \textit{subriemannian metric}. An absolutely continuous curve $\gamma:[0,1]\to M$ on the manifold $M$ is called \textit{horizontal} if it is almost everywhere tangent to the distribution $\Delta$. Using the inner product $\left<\cdot,\cdot\right>$, we can define the length $l(\gamma)$ of a horizontal curve $\gamma$ by
\[
l(\gamma)=\int_0^1|\dot\gamma(t)|dt,
\]
where $|\cdot|$ denotes the norm of the subriemannian metric $\left<\cdot,\cdot\right>$.

The subriemannian or Carnot-Caratheodory distance $d$ between two points $x$ and $y$ on the manifold $M$ is defined by
\begin{equation}\label{CCdistance}
 d(x,y)=\inf l(\gamma),
\end{equation}
where the infimum is taken over all horizontal curves which start from $x$ and end at $y$.

The above distance function may not be well-defined since there may exist two points which are not connected by any horizontal curve. For this we assume that the distribution $\Delta$ is bracket-generating. Before defining what a bracket-generating distribution is, let us introduce several notions. Let $\Delta_1$ and $\Delta_2$ be two distributions on a manifold $M$, and let $\mathfrak X(\Delta_i)$ be the space of all vector fields contained in the distribution $\Delta_i$. The distribution formed by the Lie brackets of the elements in $\mathfrak X(\Delta_1)$ with those in $\mathfrak X(\Delta_2)$ is denoted by $[\Delta_1,\Delta_2]$. More precisely,
\[
[\Delta_1,\Delta_2]_x=\textbf{span}\{w_1(x),[w_2,w_3](x)|w_i\in\mathfrak X(\Delta_j),i=1,2,3, j=1,2\}.
\]
We define inductively the following distributions: $[\Delta,\Delta]=\Delta^2$ and $\Delta^k=[\Delta,\Delta^{k-1}]$. A
distribution $\Delta$ is called \textit{$k$-generating} if $\Delta^k=TM$
and the smallest such $k$ is called the \textit{degree of nonholonomy}.
Finally the distribution is called \textit{bracket-generating} if it is
$k$-generating for some $k$.

Under the bracket-generating assumption, the subriemannian distance is well-defined thanks to the following famous Chow-Rashevskii Theorem (see \cite[Chapter 2]{Mo} for a proof):

\begin{thm}(Chow-Rashevskii)
 Assume that the manifold $M$ is connected and the distribution $\Delta$ is bracket-generating. Then there is a horizontal curve joining any two given points.
\end{thm}

Finally, let us discuss the subriemannian geodesics and the corresponding geodesic flow. As in Riemannian geometry, horizontal curves which realize the infimum in (\ref{CCdistance}) are called length minimizing geodesics (or simply geodesics). From now on, all subriemannian manifolds are assumed to be complete as a metric space. It follows that given any two points on the manifold, there is at least one constant speed geodesic joining them. Next we will discuss one type of geodesics called normal geodesics. For this let us recall several notions in the symplectic geometry of the cotangent bundle $T^*M$. Let $\pi:T^*M\to M$ be the projection map, the tautological one-form $\theta$ on $T^*M$ is defined by
\[
 \theta_\alpha(X)=\alpha(d\pi(X)),
\]
where $\alpha$ is in the cotangent bundle $T^*M$ and $X$ is a tangent vector on the manifold $T^*M$ at $\alpha$.

The symplectic two-form $\omega$ on $T^*M$ is defined as the exterior derivative of the tautological one-form: $\omega=d\theta$. It is nondegenerate in the sense that $\omega(X,\cdot)=0$ if and only if $X=0$. Given a function $H:T^*M\to\Real$ on the cotangent bundle, the Hamiltonian vector field $\vec H$ is defined by $i_{\vec H}\omega=-dH$. By the nondegeneracy of the symplectic form $\omega$, the Hamiltonian vector field $\vec H$ is uniquely defined.

Given a distribution $\Delta$ and a subriemannian metric $\left<\cdot,\cdot\right>$ on it, we can associate with it a Hamiltonian $H$, called subriemannian Hamiltonian, on the cotangent bundle $T^*M$. To do this, let $\alpha$ be in the cotangent space $T^*_xM$ at the point $x$. The subriemannian metric $\left<\cdot,\cdot\right>$ defines a bundle isomorphism $I:\Delta^*\to\Delta$ between the distribution $\Delta$ and its dual $\Delta^*$. It is defined by
\[
\left<I(\beta),\cdot\right>=\beta(\cdot),
\]
where $\beta$ is an element in the dual bundle $\Delta^*$ of the distribution $\Delta$.

By restricting the domain of the covector $\alpha$ to the subspace $\Delta_x$ of the tangent space $T_xM$, it defines an element, still called $\alpha$, in the dual space $\Delta^*$. Therefore, $I(\alpha)$ is a tangent vector contained in the space $\Delta_x$ and the subriemannian Hamiltonian $H$ corresponding to the subriemannian metric $\left<\cdot,\cdot\right>$ is defined by
\[
 H(\alpha):=\frac{1}{2}\alpha(I(\alpha))=\frac{1}{2}\left<I(\alpha),I(\alpha)\right>.
\]
Note that this construction defines the usual kinetic energy Hamiltonian in the Riemannian case.

Let $\vec H$ be the Hamiltonian vector field corresponding to the subriemannian Hamiltonian $H$ and we denote the corresponding flow, the subriemannian geodesic flow, by $e^{t\vec H}$. If $t\mapsto e^{t\vec H}(\alpha)$ is a trajectory of the subriemannian geodesic flow, then its projection $t\mapsto \gamma(t)=\pi(e^{t\vec H}(\alpha))$ is a locally minimizing geodesic. That means sufficiently short segment of the curve $\gamma$ is a minimizing geodesic between its endpoints. The minimizing geodesics obtained this way are called normal geodesics. In the special case where the distribution $\Delta$ is the whole tangent bundle $TM$, the distance function (\ref{CCdistance}) is the usual Riemannian distance and all geodesics are normal. The same is true for subriemannian manifolds, called contact subriemannian manifolds (see Section \ref{contact} for the definition), studied in this paper. However, this is not the case for general subriemannian manifolds. To introduce another class of geodesics, consider the space $\Omega$ of horizontal curves with square integrable derivatives. The endpoint map $end:\Omega\to M$ is defined by taking an element $\gamma$ in space of curves $\Omega$ and giving the endpoint $\gamma(1)$ of the curve: $end(\gamma)=\gamma(1)$. Geodesics which are regular points of the endpoint map are automatically normal and those which are critical points are called abnormal. However, there are geodesics which are both normal and abnormal (see \cite[Chapter 3]{Mo} and reference therein for more detail about abnormal geodesics).

\smallskip

\section{Contact Subriemannian and Sasakian Manifolds}\label{contact}

In this section, we recall the definition of contact subriemannian manifolds which is the main object of study for this paper. We will also recall the definition of Sasakian manifolds which served as key examples of various results. Finally we will mention some explicit examples in the three dimensional case.

A distribution $\Delta$ on a manifold $M$ is \textit{contact} if there exists a 1-form $\alpha_0$, called \textit{contact form}, for which
\begin{itemize}
\item the kernel of $\alpha_0$ is $\Delta$ (i.e. $\alpha_0(v)=0$ for each $v$ in $\Delta$) and
\item the differential $d\alpha_0$ is non-degenerate on $\Delta$ (i.e. $d\alpha_0(v,\cdot)\equiv 0$ if and only if $v=0$).
\end{itemize}
Note that the second condition implies the manifold $M$ is odd dimensional. Therefore, we can assume that the dimension of the manifold is $2n+1$. Once a contact distribution is fixed, there are lots of contact form associated with it. However, if a subriemannian metric $\left<\cdot,\cdot\right>$ is also fixed on the distribution, then there is a unique contact form $\alpha_0$ such that the restriction of the $2n$-form $d\alpha_0\wedge ... \wedge d\alpha_0$ to the distribution $\Delta$ coincides with the volume form induced by the subriemannian metric $\left<\cdot,\cdot\right>$ on $\Delta$. Therefore, we say that the subriemannian manifold $(M,\Delta,\left<\cdot,\cdot\right>)$ is a contact subriemannian manifold if $\Delta$ is a contact distribution and we call the 1-form $\alpha_0$ defined above the \textit{induced contact form} of $(M,\Delta,\left<\cdot,\cdot\right>)$.

For each contact subriemannian manifold $(M,\Delta,\left<\cdot,\cdot\right>)$, we can associate with it a unique vector field $v_0$, called the \textit{Reeb field}. If $\alpha_0$ is the induced contact form, then $v_0$ is defined by conditions $\alpha_0(v_0)=1$ and $d\alpha_0(v_0,\cdot)=0$. Note that the first condition implies the Reeb field $v_0$ is transversal to the distribution $\Delta$.

Using the Reeb field $v_0$, we can define a natural measure on the subriemannian manifold. Let $v_1,...,v_{2n}$ be a basis in the contact distribution $\Delta$ which is orthonormal with respect to the given subriemannian metric. Let $\eta$ be the $(2n+1)$-form defined by the condition $\eta(v_0,...,v_{2n})=1$. The measure induced by this volume form $\eta$, which will be denoted by the same symbol throughout this paper, is an example of a Popp's measure. Popp's measures can be defined for any subriemannian manifold. For the detail definition of this measure in general, see \cite[Chapter 10]{Mo}. From now on, when we consider a contact subriemannian manifold as a metric measure space, it always refers to the triple $(M,d,\eta)$ where $d$ is the subriemannian distance and $\eta$ is the Popp's measure.

Before giving examples of contact subriemannian manifolds, let us recall the definition of an important class of manifolds, called \textit{Sasakian manifolds}. A three dimensional contact manifold is Sasakian if the Reeb field $v_0$ is a subriemannian isometry. Note that by the definition of the Reeb field $v_0$, the flow of $v_0$ preserves the induced contact form and hence the distribution $\Delta$. Therefore, subriemannian isometry here means the flow of the Reeb field $v_0$ preserves the subriemannian length of tangent vectors in $\Delta$. Higher dimensional contact manifolds for which the Reeb field $v_0$ is a subriemannian isometry are called K-contact. Sasakian manifolds are defined by the integrability of certain tensor. In general, Sasakian manifolds are K-contact, but not conversely. In the three dimensional case, the two notions coincide (see \cite{Bl} for the detail).

If we assume further that the Reeb field $v_0$ of the Sasakian manifold generates a free and proper group action (i.e. the flow of $v_0$ is a free and proper group action), then the quotient $N:=M/G$ of the manifold $M$ by this $G$-action ($G=S^1\text{ or } \Real$) is again a manifold. Let $\pi_M:M\to N$ be the quotient map. Then there is a Riemannian metric on $N$ such that the restriction of $d\pi$ to the distribution $\Delta$ is an isometry. Many of the interesting examples of subriemannian manifolds have this structure. The Heisenberg group $\mathbb H^n$ is one of them.

The standard subriemannian structure of the Heisenberg group $\mathbb H^n$ can be defined as follows. The underlying manifold of $\mathbb H^n$ is the $2n+1$-dimensional Euclidean space $M=\Real^{2n+1}$. If we denote the coordinates of this Euclidean space by $\{x_0,x_1,...,x_{2n}\}$, then the distribution $\Delta$ is defined by
\[
\Delta=\textbf{span}\left\{X_i,Y_i|i=1,...,n\right\},
\]
where $X_i=\partial_{x_i}-\frac{1}{2}x_{n+i}\partial_{x_0}$ and $Y_i=\partial_{x_{n+i}}+\frac{1}{2}x_{i}\partial_{x_0}$.

The standard subriemannian metric $\left<\cdot,\cdot\right>$ is the one for which the vector fields $\left\{X_i,Y_i|i=1,...,n\right\}$ are orthonormal. The induced contact form $\alpha_0$ is given by
\[
\alpha_0=-dx_0+\frac{1}{2}\sum_{i=1}^n(x_idx_{n+i}-x_{n+i}dx_i).
\]
The Reeb field $v_0$ in this case is $-\partial_{x_0}$ and the Popp's measure $\eta$ is the $2n+1$-dimensional Lebesgue measure. The Reeb field $v_0$, in this case, is a subriemannian isometry. It also defines a proper $\Real$-action and the quotient manifold $N$ is the $2n$-dimensional Euclidean space $\Real^{2n}$. The standard subriemannian structure $\left<\cdot,\cdot\right>$ on $\mathbb H^n$ descends to the standard Euclidean structure on $\Real^{2n}$.

We end this section with two more examples of contact subriemannian manifolds with the above symmetry structure. For more examples of subriemannian manifolds with symmetry, see \cite[Chapter 11]{Mo}. For other examples of contact manifolds, see \cite{Bl}.

Recall that $SU(2)$, the special unitary group, consists of $2\times2$ unitary matrices. The Lie algebra $su(2)$ consists of skew Hermitian matrices with trace zero. The left invariant vector fields of the following two elements in $su(2)$
\[
v_1=\left(
\begin{array}{cc}
0 & 1/2\\
-1/2 & 0
\end{array}\right),\quad v_2=\left(
\begin{array}{cc}
0 & i/2\\
i/2 & 0
\end{array}\right)
\]
span the standard distribution $\Delta$ on $SU(2)$. The standard subriemannian metric is given by the condition $\left<v_i,v_j\right>=\delta_{ij}$, $i=1,2$. The Reeb field $v_0$ is given by
\[
v_0=\left(
\begin{array}{cc}
-i/2 & 0\\
0 & i/2
\end{array}\right).
\]

The flow of the Reeb field defines a $S^1$-action on $SU(2)$ called the Hopf fibration. The quotient $N$ of $SU(2)$ by this action is the standard 2-sphere $S^2$. The standard subriemannian metric on $SU(2)$ descends to the Riemannian metric on $S^2$ of contant curvature 1.

The special linear group $SL(2)$ is the set of all $2\times 2$ matrices with real coefficients and determinant 1. The Lie algebra $sl(2)$ is the set of all $2\times 2$ real matrices with trace zero. The left invariant vector fields of the following two elements in $sl(2)$
\[
v_1=\left(
\begin{array}{cc}
1/2 & 0\\
0 & -1/2
\end{array}\right),\quad v_2=\left(
\begin{array}{cc}
0 & 1/2\\
1/2 & 0
\end{array}\right)
\]
span the standard distribution $\Delta$ on $SL(2)$. The standard subriemannian metric on $SL(2)$ is defined by $\left<v_i,v_j\right>=\delta_{ij}$, $i=1,2$. The Reeb field in this case is $v_0$, where
\[
v_0=\left(
\begin{array}{cc}
0 & -1/2\\
1/2 & 0
\end{array}\right).
\]

The flow of the Reeb field also defines a $S^1$-action on $SU(2)$. The quotient $N$ of $SL(2)$ by this action is the upper half-space with the standard non-Euclidean structure.

\smallskip

\section{Introduction to Optimal Transportation}\label{OptimalIntro}

In this section, we give a quick introduction to the optimal transportation problem. A standard reference on this is the book \cite{Vi2}.

Let $M$ be a metric space with distance function $d$. Let $\mu_0$ and $\mu_1$ be two Borel probability measures on $M$. The theory of optimal transportation starts with the following minimization problem
\begin{equation}\label{optimalorigin}
\inf_{\varphi_*\mu_0=\mu_1}\int_Md^2(x,\varphi(x))\, d\mu_0(x)
\end{equation}
where the infimum is taken over all Borel maps $\varphi:M\to M$ which pushes $\mu_0$ forward to $\mu_1$ (i.e. $\mu_0(\varphi^{-1}(U))=\mu_1(U)$ for all Borel sets $U$ in $M$).

By the famous work of \cite{Ka}, the relaxed version of the above problem given below in (\ref{optimal}) always has a solution (i.e. existence of minimizer).

\begin{equation}\label{optimal}
\inf_{(\pi_1)_*\Pi=\mu_0, (\pi_2)_*\Pi=\mu_1}\int_{M\times M}d^2(x,y)\, d\Pi(x,y)
\end{equation}
where $\pi_1,\pi_2:M\times M\to M$ are projections onto the first and second component, respectively, and
the infimum is taken over all Borel measures $\Pi$ on $M\times M$ satisfying $(\pi_1)_*\Pi=\mu_0$ and $(\pi_2)_*\Pi=\mu_1$ (i.e. $\Pi(U\times M)=\mu_0(U)$ and $\Pi(M\times U)=\mu_1(U)$ for all Borel sets $U$ in $M$).

The existence and uniqueness of solution to the origin problem (\ref{optimalorigin}) were proved much later in (\cite{Br}) in the Euclidean setting under certain assumptions on the measures $\mu_0$ and $\mu_1$. It was later extended to the compact Riemannian setting by \cite{Mc2}. The following is a summary of their results (see also the result in \cite{Fi} where all the compactness assumptions are removed).

\begin{thm}\cite{Br,Mc2}\label{BrMc}
Let $M$ be a Riemannian manifold with Riemannian distance $d$. Assume that the measures $\mu_0$ and $\mu_1$ have compact supports and the measure $\mu_0$ is absolutely continuous with respect to the Riemannian volume. Then the optimal transportation problem (\ref{optimalorigin}) has a solution $\varphi$ which is unique up to a set of $\mu$-measure zero. Moreover, there exists a Lipschitz function $\mathfrak f:M\to\Real$ such that the map $\varphi$ is given by
\[
\varphi(x)=\exp(\nabla \mathfrak f(x)).
\]
\end{thm}

The problem (\ref{optimalorigin}) in the subriemannian setting was first considered in \cite{AmRi}. Under the same  assumptions as in Theorem \ref{BrMc} on the measures, the existence and uniqueness of the solution was shown when the space is the Heisenberg group equipped with the standard subriemannian metric (see Section \ref{subgeo} for the definition). The generalization to more general subriemannian manifolds is later done in \cite{AgLe1}. In \cite{AgLe1}, the authors proved that the existence and uniqueness theorem holds when the subriemannian manifold is 2-generating (see Section \ref{subgeo} for the definition of $k$-generating distribution). In particular, it is applicable to the contact subriemannian case considered in the present work.

\begin{thm}\cite{AgLe1}\label{subtransport}
Let $M$ is a subriemannian manifold with subriemannian distance $d$ and a 2-generating distribution $\Delta$. Assume that the measures $\mu_0$ and $\mu_1$ have compact supports and the measure $\mu_0$ is absolutely continuous with respect to a Riemannian volume. Then the optimal transportation problem (\ref{optimalorigin}) has a solution $\varphi$ which is unique up to a set of $\mu_0$-measure zero. Moreover, there exists a function $\mathfrak f:M\to\Real$ which is Lipschitz with respect to a Riemannian metric such that the map $\varphi$ is given by
\[
\varphi(x)=\pi(e^{1\cdot\vec H} (d\mathfrak f_x)).
\]
where $e^{t\vec H}$ denotes the subriemannian geodesic flow and $\pi:T^*M\to M$ is the natural projection.
\end{thm}

The difficulty in extending the above theorem to all subriemannian manifolds lies in the presence of abnormal minimizers. Using geometric measure theory, \cite{FiRi} is able to extend Theorem \ref{subtransport} to more general subriemannian manifolds. However, the problem of showing uniqueness or non-uniqueness of solutions to (\ref{optimalorigin}) in the subriemannian case remains unsolved in general.

\smallskip

\section{Optimal Transportation and Ricci Curvature}\label{optimalRicci}

In this section, we give a very brief overview of results concerning the connection of optimal transportation with  generalized Ricci curvature lower bound (see \cite{LoVi1,LoVi2,St1,St2} for a detail discussion).

The optimal transportation problem in (\ref{optimal}) defines  a distance function on the space of all Borel probability measures of a given metric space. More precisely, let $\mathfrak X$ be a locally compact complete separable metric space with distance function $\mathfrak d$. Let $\mathcal P$, called the Wasserstein space, be the space of all Borel probability measures $\mu$ of $\mathfrak X$ such that the following integral is finite for some point $x_0$ in $\mathfrak X$
\[
\int_{\mathfrak X}\mathfrak d^2(x,x_0)d\mu(x).
\]

The Wasserstein distance function $\mathcal W$ on $\mathcal P$ is defined by the optimal transportation problem as follows.
\begin{equation}\label{WassersteinD}
\mathcal W(\mu_0,\mu_1)=\left(\inf_{(\pi_1)_*\Pi=\mu_0, (\pi_2)_*\Pi=\mu_1}\int_{\mathfrak X\times \mathfrak X}\mathfrak d^2(x,y)\, d\Pi(x,y)\right)^{1/2}.
\end{equation}

Assume that the space $\mathfrak X$ is a geodesic space. Then the Wasserstein space $\mathcal P$ equipped with the Wasserstein distance $\mathcal W$ is a geodesic space (i.e. distance between two points is given by the length of the shortest curve, called geodesic, connecting them, see \cite{St2} for the precise definition of geodesic space and the proof of this fact).

\begin{rem}\label{displace}
Assume that the metric space $(\mathfrak X,\mathfrak d)$ is a contact subriemannian manifold and the measure $\mu_0$ is absolutely continuous with respect to a Riemannian volume. Then the geodesics of the corresponding Wasserstein distance are given by
\[
t\mapsto (\varphi_t)_*\mu_0
\]
where $\varphi_t=\pi(e^{t\vec H}(d\mathfrak f))$ and $\mathfrak f$ is defined as in Theorem \ref{subtransport}. These paths of measures, called displacement interpolations, were first introduced in \cite{Mc1}.
\end{rem}

Finally let us fix a locally finite measure $\nu$ and introduce the relative entropy functional $\ent\mu\nu$ on $\mathcal P$ by
\[
\ent \mu\nu=
\begin{cases}
\int_M g\log g \, d\nu & \text{if } \mu=g\nu\\
+\infty & \text{otherwise}.
\end{cases}
\]

Formally, a metric measure space $(\mathfrak X,\mathfrak d,\nu)$ satisfies the curvature-dimension condition $CD(K,\infty)$ if the above relative entropy functional has second derivative bounded below by $K$ along any geodesic in the Wasserstein space $\mathcal P$ of $M$ equipped with the Wasserstein distance $\mathcal W$. More precisely, the second derivative is replaced by the following difference quotient.

\begin{defn}
The metric measure space $(\mathfrak X,\mathfrak d,\nu)$ satisfies \textit{curvature-dimension condition} $CD(K,\infty)$ if, given any two measures $\mu_0$ and $\mu_1$, there is a geodesic $\mu_t$ in $\mathcal P$ such that the followings hold for all $t$ in $[0,1]$
\[
\frac{K}{2}t(1-t)\mathcal W^2(\mu_0,\mu_1) \leq (1-t)\ent{\mu_0}\nu+t\ent{\mu_1}\nu -\ent{\mu_t}\nu.
\]
\end{defn}

There are also other curvature-dimension conditions $CD(K,N)$ for $N>0$ finite. The definitions of these conditions are similar to that of $CD(K,\infty)$ but are more involved. Their detail definitions as well as related results can found in \cite{St2}.

In the Riemannian case, the condition $CD(K,\infty)$ is the same as Ricci curvature bounded below by $K$. More precisely, the following holds.

\begin{thm}\cite{OtVi,CoMcSc,StVo}
Assume that $\mathfrak X$ is a complete Riemannian manifold equipped with the Riemannian distance $\mathfrak d$ and the measure $\nu$ induced by the Riemannian volume form. If we denote the Ricci curvature by $\textbf{Ric}$ and the Riemannian metric by $\left<\cdot,\cdot\right>$, then the metric measure space $(\mathfrak X,\mathfrak d,\nu)$ satisfies curvature-dimension condition $CD(K,\infty)$ if and only if
\[
\textbf{Ric}(v,v)\geq K|v|^2
\]
for all tangent vector $v$ in the tangent bundle $T\mathfrak X$.
\end{thm}

Besides Riemannian manifolds, it was shown in \cite{Oh2} that $CD(K,\infty)$ is equivalent to the flag Ricci tensor bounded below by $K$ in the Finsler case. The situation in the subriemannian case is completely different. In \cite{Ju}, it was shown that the most basic subriemannian example, the Heisenberg group (see Section \ref{contact} for the precise definition of the Heisenberg group and the standard subriemannian structure on it), does not satisfy any curvature-dimension condition mentioned above (however see \cite{BaGa1,BaBoGa,BaGa2} for a different curvature-dimension condition in the subriemannian setting which is satisfied by the Heisenberg group in particular). On the other hand, it was shown in \cite{Ju} that the Heisenberg group satisfies the measure contraction property $MCP(K,N)$ defined below. The following is the definition of measure contraction property (see \cite{LoVi1,Oh1,St1} for more details).

\begin{defn}\label{tinter}
Let $\mathfrak X$ be a geodesic space. A point $z$ is a $t$-intermediate point of $x$ and $y$ if there is a geodesic $\gamma:[0,1]\to\mathfrak X$ such that $\gamma(0)=x$, $\gamma(1)=y$, and $\gamma(t)=z$.
\end{defn}

\begin{defn}\label{MCP}
The metric measure space $(\mathfrak X,\mathfrak d,\nu)$ satisfies the measure contraction property $MCP(K,N)$ if for each $t$ in $(0,1)$, there is a Markov kernel $P_t$ which takes a point in $\mathfrak X\times\mathfrak X$ to a measure in $\mathfrak X$ such that for $\nu^2$ almost every $(x_0, x)$ and for $P_t(x_0, x)$ almost every $z$ the point $z$ is a $t$-intermediate point of $x_0$ and $x$, and the followings hold:
\[
\nu(U)\geq  \int_{\mathfrak X}(1-t) \left(\frac{s_K((1-t)D(x))}{s_K(D(x))}\right)^{N-1}P_t(x,x_0)(U) d\nu(x)
\]
and
\[
\nu(U)\geq  \int_{\mathfrak X}t \left(\frac{s_K(tD(x))}{s_K(D(x))}\right)^{N-1}P_t(x_0,x)(U) d\nu(x)
\]
for any measurable set $U$ and $\nu$-almost every $x$ in $M$, where
\[
s_K(r)=
\begin{cases}
\frac{1}{\sqrt K}\sin(\sqrt K \,r) & \text{if }K>0\\
r & \text{if }K=0\\
\frac{1}{\sqrt {-K}}\sinh(\sqrt {-K}\, r) & \text{if }K<0
\end{cases}
\]
and $D(x)=\frac{d(x_0,x)}{\sqrt{N-1}}$.
\end{defn}

As mentioned in the introduction, $MCP(K,N)$ is another characterization of Ricci curvature lower bound for $N$-dimensional Riemannian manifolds.

\begin{thm}\cite{St2,Oh1}
Assume that $\mathfrak X$ is a $N$-dimensional complete Riemannian manifold equipped with the Riemannian distance $\mathfrak d$ and the measure $\nu$ induced by the Riemannian volume form. If we denote the Ricci curvature by $\textbf{Ric}$ and the Riemannian metric by $\left<\cdot,\cdot\right>$, then the metric measure space $(\mathfrak X,\mathfrak d,\nu)$ satisfies the measure contraction property $MCP(K,N)$ if and only if
\[
\textbf{Ric}(v,v)\geq K|v|^2
\]
for all tangent vector $v$ in the tangent bundle $T\mathfrak X$.
\end{thm}

Finally we end this section with the following theorem, proved in \cite{Ju}, which motivates the present work.

\begin{thm}\cite{Ju}
Let $\mathbb H^n$ be the $2n+1$ dimensional Heisenberg group. Let $\mathfrak d$ be the standard subriemannian distance and $\nu$ be the $(2n+1)$-dimensional Lebesgue measure $dx^{2n+1}$. Then the metric measure space $(\mathbb H^n,\mathfrak d,dx^{2n+1})$ satisfies the measure contraction property $MCP(0,2n+3)$.
\end{thm}

\smallskip

\section{Generalized Curvatures on Subriemannian Manifolds}\label{gencur}

In this section, we recall the definition of the curvature type invariants studied in \cite{AgGa,AgZe,LiZe1,LiZe2} and specialize it to the case of a three dimensional contact subriemannian manifold.

Let $e^{t\vec H}$ be the subriemannian geodesic flow defined in Section \ref{subgeo} and let $\alpha$ be a point in the manifold $T^*M$. As mentioned in the introduction, the idea is to construct a Fr\'enet-Serret type frame along the curve $t\mapsto e^{t\vec H}(\alpha)$ so that the pulled back frame, called canonical Darboux frame, satisfies certain differential equations, called structural equations. The coefficients of these equations, in turn, defines the curvature operator that we need.

The vertical space $V_\alpha$ at $\alpha$ of the bundle $\pi:T^*M\to M$ is defined as the kernel of the map $d\pi_\alpha:T_\alpha T^*M\to T_{\pi(\alpha)}M$. Recall that a subspace $V$ of a symplectic vector space of dimension $2m$ is Lagrangian if the symplectic form restricted to $V$ vanishes and the dimension of $V$ is $m$. Each of these vertical spaces $V_\alpha$ is a Lagrangian subspace with respect to the canonical symplectic form $\omega$ defined in Section \ref{subgeo}. On the other hand, the differential $de^{-t\vec H}:T_{e^{t\vec H}(\alpha)} T^*M\to T_\alpha T^*M$ of the map $e^{-t\vec H}$ is a symplectic transformation (i.e. it preserves the symplectic form) between the symplectic vector spaces $T_{e^{t\vec H}(\alpha)} T^*M$ and $T_\alpha T^*M$. Therefore, the one parameter family of subspaces
\[
t\mapsto J_\alpha(t):=de^{-t\vec H}(V_{e^{t\vec H}(\alpha)})
\]
defines a curve of Lagrangian subspaces contained in a single symplectic vector space $T_\alpha T^*M$. This curve is called the \textit{Jacobi curve} at $\alpha$.

Recall that the space of all Lagrangian subspaces in a symplectic vector space $\Sigma$ is a finite dimensional manifold (in fact a homogeneous space of the symplectic group), called the Lagrangian Grassmannian $LG(\Sigma)$ of $\Sigma$. The Jacobi curve defined above is a smooth curve in the Lagrangian Grassmannian $LG(T_\alpha T^*M)$. The curvature type invariants of the geodesic flow $e^{t\vec H}$ are simply differential invariants of the Jacobi curve under the action of the symplectic group (see \cite{LiZe1,LiZe2} for further details). The construction of differential invariants for a general curve $t\mapsto J(t)$ in the Lagrangian Grassmannian $LG(\Sigma)$ of a symplectic vector space $\Sigma$ was done in the recent papers \cite{LiZe1,LiZe2}, though partial results were obtained earlier (see \cite{Ah,Fo,Gr,AgGa,AgZe}). 

Recall that a basis $\{e_1,...,e_n,f_1,...,f_n\}$ in a symplectic vector space with a symplectic form $\omega$ is a Darboux basis if it satisfies $\omega(e_i,e_j)=\omega(f_i,f_j)=0$, and $\omega(f_i,e_j)=\delta_{ij}$. Given a subriemannian Hamiltonian, there is a moving Darboux basis $\{e_1(t),...,e_n(t),f_1(t),...,f_n(t)\}$, called canonical Darboux frame, of the symplectic vector space $T_\alpha T^*M$ such that $J_\alpha(t)=\textbf{span}\{e_1(t),...,e_n(t)\}$ and, more importantly, the canonical Darboux frame satisfies a system of first order ODEs of specific form, called structural equations. This defines a splitting of the symplectic vector space $T_\alpha T^*M=J_\alpha(t)\oplus\hat J_\alpha(t)$, where $\hat J_\alpha(t)=\textbf{span}\{f_1(t),...,f_n(t)\}$. In particular, the subspace $J_\alpha(0)$ is the vertical space $V_\alpha$ of the bundle $\pi:T^*M\to M$ and the subspace $\hat J_\alpha(0)$ is a complimentary
subspace to $J_\alpha(0)=V_\alpha$ at time $t=0$. Hence, $\bigcup_{\alpha\in T^*M}\hat J_\alpha(0)$ defines an Ehresmann connection on the bundle $\pi:T^*M\to M$.

In the Riemannian case, this is, under the identification of the tangent and cotangent spaces by the Riemannian metric, simply the Levi-Civita connection (see \cite[Proposition 5.2]{AgGa}). The canonical Darboux frame, in this case, satisfies the following equations which is the Jacobi field equation (up to certain identifications of tangent and cotangent spaces)
\[
\dot e_i(t)=f_i(t),\quad \dot f_i(t)=-R^{ij}_\alpha(t)e_j(t).
\]
The matrix $R_\alpha:=R_\alpha(0)$ with $ij$-th entries given by $R^{ij}_\alpha(0)$ above is the Riemannian curvature operator (again up to certain identifications).

Using the above splitting we can also define a generalization of the Ricci curvature in the Riemannian geometry. Indeed let $\pi_{J_\alpha(t)}$ and $\pi_{\hat J_\alpha(t)}$ be the projections, corresponding to the splitting $T_\alpha T^*M=J_\alpha(t)\oplus\hat J_\alpha(t)$, onto the subspaces $J_\alpha(t)$ and $\hat J_\alpha(t)$, respectively. Let $w(\cdot)$ be a path contained in the Jacobi curve $J_\alpha(\cdot)$ (i.e. $w(t)\in J_\alpha(t)$ for all $t$). Then the projection $\pi_{\hat J_\alpha(t)}\dot w(t)$ of its derivative $\dot w(t)$ onto the subspace $\hat J_\alpha(t)$ depends only on the vector $w(t)$ but not on the curve $w(\cdot)$. Therefore, it defines a linear operator $\Phi^t_{J_\alpha\hat J_\alpha}:J_\alpha(t)\to\hat J_\alpha(t)$
\[
 \Phi^t_{J_\alpha\hat J_\alpha}(w(t))=\pi_{\hat J_\alpha(t)}\left(\dot w(t)\right).
\]
Similarly we can also define another operator $\Phi_{\hat J_\alpha J_\alpha}^t:\hat J_\alpha(t)\to J_\alpha(t)$ by switching the role of $J$ and $\hat J$ above. The composition of $\Phi_{\hat J_\alpha J_\alpha}^0$ and $\Phi_{ J_\alpha\hat J_\alpha}^0$ defines a linear operator $\Phi_{\hat J_\alpha J_\alpha}^0\circ\Phi_{ J_\alpha\hat J_\alpha}^0:J_\alpha(0)=V_\alpha\to V_\alpha$ of the vertical space $V_\alpha$. Finally the \textit{generalized Ricci curvature} $\ric(\alpha)$ at $\alpha$ is defined by the negative of the trace of $\Phi_{\hat J_\alpha J_\alpha}^0\circ\Phi_{ J_\alpha\hat J_\alpha}^0$. When the geodesic flow $e^{t\vec H}$ is Riemannian, the generalized Ricci curvature $\ric$ reduces to the usual Ricci curvature (under certain identifications of tangent and cotangent spaces).

Now let us consider the three dimensional contact subriemannian case. The structural equations, in this case, have the following form (see Section \ref{proofstructural} for the proof):

\begin{thm}\label{structural}
Let $(M,\Delta,\left<\cdot,\cdot\right>)$ be a three dimensional contact subriemannian manifold. For each fixed $\alpha$ in $T^*M$, there is a moving Darboux frame
\[
e_1(t),e_2(t),e_3(t),f_1(t),f_2(t),f_3(t)
\]
of the symplectic vector space $T_\alpha T^*M$ and functions $R^{11}_\alpha(t),R^{22}_\alpha(t)$ of time $t$ such that $\{e_1(t),e_2(t),e_3(t)\}$ form a basis for the Jacobi curve $J_\alpha(t)$ and it satisfies the following structural equations
\[
\left\{
  \begin{array}{ll}
    \dot e_1(t)=f_1(t), \\
    \dot e_2(t)=e_1(t), \\
    \dot e_3(t)=f_3(t), \\
    \dot f_1(t)=-R^{11}_\alpha(t)e_1(t)-f_2(t),\\
    \dot f_2(t)=-R^{22}_\alpha(t)e_2(t),\\
    \dot f_3(t)=0.
  \end{array}
\right.
\]

Moreover, the generalized Ricci curvature $\mathfrak {Ric}(\alpha)$ at $\alpha$ is given by $\mathfrak {Ric}(\alpha)=R^{11}_\alpha(0)$.
\end{thm}

Next we will write down explicit formulas (Theorem \ref{structuraldetail}) for the canonical Darboux frame and the differential invariants $R^{11}(t)$ and $R^{22}(t)$ in Theorem \ref{structural}. Let $\{v_1,v_2\}$ be a local orthonormal frame in the contact distribution $\Delta$ with respect to the subriemannian metric $\left<\cdot,\cdot\right>$ and let $v_0$ be the Reeb field. This defines a convenient frame $\{v_0,v_1,v_2\}$ in (a neighborhood of) the tangent bundle $TM$ and we let  $\{\alpha_0,\alpha_1,\alpha_2\}$ be the corresponding dual co-frame in the cotangent bundle $T^*M$ (i.e. $\alpha_i(v_j)=\delta_{ij}$).

The frame $\{v_0,v_1,v_2\}$ and the co-frame $\{\alpha_0,\alpha_1,\alpha_2\}$ defined above induces a frame in the tangent bundle $TT^*M$ of the cotangent bundle $T^*M$. Indeed, let $\vec\alpha_i$ be the vector fields on the cotangent bundle $T^*M$ defined by $i_{\vec\alpha_i}\omega=-\alpha_i$. Note that the symbol $\alpha_i$ in the definition of $\vec\alpha_i$ represents the pull back $\pi^*\alpha_i$ of the 1-form $\alpha$ on the manifold $M$ by the projection $\pi:T^*M\to M$. This convention of identifying forms in the manifold $M$ and its pull back on the cotangent bundle $T^*M$ will be used for the rest of this paper without mentioning. Note that we use the same symbol $\alpha$ to represent a 1-form on $M$ and also its pull back $\pi^*\alpha$ on $T^*M$. It will be clear from the context which geometric object $\alpha$ represents.

Let $h_i:T^*M\to\Real$ be the Hamiltonian lift of the vector fields $v_i$, defined by $h_i(\alpha)=\alpha(v_i)$. Let $\vec\xi_1$ and $\vec\xi_2$ be the vector fields defined by $\vec\xi_1=h_1\vec\alpha_2-h_2\vec\alpha_1$ and $\vec\xi_2=h_1\vec\alpha_1+h_2\vec\alpha_2$. Then the vector fields $\vec h_0, \vec h_1,\vec h_2, \vec\alpha_0, \vec\xi_1, \vec\xi_2$ define a local frame for the tangent bundle $TT^*M$ of the cotangent bundle $T^*M$. We are going to write the canonical Darboux frame in terms of this convenient local frame. Finally, we also let $h_{ij}$ be the Hamiltonian lift of $[v_i,v_j]$ defined by $h_{ij}(\alpha)=\alpha([v_i,v_j])$.

Under the above notation the subriemannian Hamiltonian is given by $H=\frac{1}{2}((h_1)^2+(h_2)^2)$ and the Hamiltonian vector field is $\vec H=h_1\vec h_1+h_2\vec h_2$. Let $\mathfrak d_s:T^*M\to T^*M$ be the dilation in the fibre direction defined by $\mathfrak d_s(\alpha)=s\alpha$ and let $\vec E$ be the Euler field defined by $\vec{E}(\alpha)=\dds\mathfrak d_s(\alpha)\Big|_{s=1}$. It is also given by $\vec E=-h_0\vec\alpha_0-\xi_2$.

We also need the bracket relations of the vector fields $v_0,v_1,v_2$. Let $c_{ij}^k$ be the functions on the manifold $M$ defined by

\begin{equation}\label{bracket}
 [v_i,v_j]= c_{ij}^0v_0+c_{ij}^1v_1+c_{ij}^2v_2.
\end{equation}

Note that $c_{ij}^k=-c_{ji}^k$. The dual version of the above relation is
\begin{equation}\label{dualbracket}
 d\alpha_k=-\sum\limits_{0\leq i<j\leq 2}c_{ij}^k\alpha_i\wedge\alpha_j.
\end{equation}

By (\ref{dualbracket}), the definition of the Reeb field $v_0$, and that of the induced from $\alpha_0$, it follows that $d\alpha_0=\alpha_1\wedge\alpha_2$. Therefore, $c_{01}^0=c_{02}^0=0$ and $c_{12}^0=-1$. If we also take the exterior derivative of the equation in (\ref{dualbracket}), we get $c_{01}^1+c_{02}^2=0$. We summarize
\begin{lem}\label{strcst}
\[
c_{01}^0=c_{02}^0=0,\quad c_{12}^0=-1,\quad c_{01}^1+c_{02}^2=0.
\]
\end{lem}

Finally we come to the main theorem of this section. Note that all vector fields in Theorem \ref{structuraldetail}, Theorem \ref{structuraldetailSas}, and their proofs should be evaluated at $\alpha$. They are omitted to avoid heavy notations.

\begin{thm}\label{structuraldetail}
The canonical Darboux frame
\[
e_1(t),e_2(t),e_3(t),f_1(t),f_2(t),f_3(t)
\]
and the differential invariants $R^{11}_\alpha(t)$ and $R^{22}_\alpha(t)$ in Theorem \ref{structural} satisfy
$R^{11}_\alpha(t)=R^{11}_{e^{t\vec H}(\alpha)}(0)$, $R^{22}_\alpha(t)=R^{22}_{e^{t\vec H}(\alpha)}(0)$, and

\[
\left\{
  \begin{array}{ll}
      e_1(t)=\frac{1}{\sqrt{2H}}(e^{t\vec H})^*\vec\xi_1,\\
    e_2(t)=\frac{1}{\sqrt{2H}}(e^{t\vec H})^*\vec\alpha_0,\\
    e_3(t)=\frac{1}{\sqrt{2H}}(e^{t\vec H})^*\vec E=\frac{1}{\sqrt{2H}}(\vec E-t\vec H),\\
    f_1(t)=\frac{1}{\sqrt{2H}}(e^{t\vec H})^*[h_1\vec h_2-h_2\vec h_1+\chi_0\vec\alpha_0+(\vec\xi_1h_{12})\vec\xi_1-h_{12}\vec\xi_2],\\
    f_2(t)=\frac{1}{\sqrt{2H}}(e^{t\vec H})^*[2H\vec h_0-h_0\vec H-\chi_1\vec\alpha_0+({\vec\xi_1} a)\vec\xi_1-a\vec\xi_2],\\
    f_3(t)=-\frac{1}{\sqrt{2H}}\vec H,\\
    \mathfrak {Ric}(\alpha):=R^{11}_\alpha(0)=h_0^2+2H\kappa-\frac{3}{2}\vec\xi_1 a,\\
    R^{22}_\alpha:=R^{22}_\alpha(0)=R^{11}_\alpha(0)\vec\xi_1 a-3\vec{H} \vec\xi_1\vec{H} a+3\vec{H}^2\vec\xi_1 a+ \vec\xi_1\vec{H}^2 a.
  \end{array}
\right.
\]
where
\[
\begin{array}{ll}
a=dh_0(\vec H),\\
\chi_0=h_2h_{01}-h_1h_{02}+\vec\xi_1a, \\
\chi_1=h_0a+2{\vec H}\vec\xi_1 a-\vec\xi_1\vec Ha,\\
\kappa=v_1c_{12}^2-v_2c_{12}^1-(c_{12}^1)^2-(c_{12}^2)^2-\frac{1}{2}(c_{01}^2-c_{02}^1),
\end{array}
\]
and $v_ic_{jk}^l$ denotes the directional derivative of the function $c_{jk}^l$ with respect to the vector field $v_i$.
\end{thm}

The proof of Theorem \ref{structuraldetail} is postponed to Section \ref{proofstructural}.

\begin{rem}
Note that $\alpha\mapsto \ric(\alpha)=R^{11}_\alpha(0)$ is a quadratic form on $T^*M$ which is positive on the kernel of the subriemannian Hamiltonian $H$. On the other hand, $\alpha\mapsto R^{22}_\alpha$ is a form of degree 4.
\end{rem}

\begin{rem}\label{TW}
It was shown in \cite{AgLe2} that $\kappa$ coincides with the Tanaka-Webster curvature in CR geometry.
\end{rem}

Recall that $a=dh_0(\vec H)$ defined in Theorem \ref{structuraldetail} is the Poisson bracket of the subriemannian Hamiltonian $H$ and the Hamiltonian lift $h_0$ of the Reeb field $v_0$. It follows immediately that a three dimensional contact subriemannian is Sasakian if and only if $a\equiv 0$. It turns out that this is also equivalent to $R^{22}\equiv 0$.

\begin{thm}\label{R22}
A three dimensional contact subriemannian manifold is Sasakian if and only if $R^{22}\equiv 0$.
\end{thm}

For the proof of this, see Section \ref{proofstructural}. In the Sasakian case, the equations in Theorem \ref{structuraldetail} simplify to

\begin{thm}\label{structuraldetailSas}
Assume that the subriemannian manifold in Theorem \ref{structural} is Sasakian. Then the canonical Darboux frame
\[
e_1(t),e_2(t),e_3(t),f_1(t),f_2(t),f_3(t)
\]
and the differential invariants $R^{11}_\alpha(t)$ and $R^{22}_\alpha(t)$ satisfy
\[
R^{11}_\alpha(t)=R^{11}_{e^{t\vec H}(\alpha)}(0),\quad R^{22}_\alpha(t)=R^{22}_{e^{t\vec H}(\alpha)}(0),
\]
\[
\mathfrak {Ric}(\alpha):=R^{11}_\alpha(0)=h_0^2+2H\kappa,\quad  R^{22}_\alpha:=R^{22}_\alpha(0)=0,
\]
and
\[
\left\{
  \begin{array}{ll}
      e_1(t)=\frac{1}{\sqrt{2H}}(e^{t\vec H})^*\vec\xi_1,\\
    e_2(t)=\frac{1}{\sqrt{2H}}(e^{t\vec H})^*\vec\alpha_0,\\
    e_3(t)=\frac{1}{\sqrt{2H}}(e^{t\vec H})^*\vec E=\frac{1}{\sqrt{2H}}(\vec E-t\vec H),\\
    f_1(t)=\frac{1}{\sqrt{2H}}(e^{t\vec H})^*[h_1\vec h_2-h_2\vec h_1+2Hc_{01}^2\vec\alpha_0+(\vec\xi_1h_{12})\vec\xi_1-h_{12}\vec\xi_2],\\
    f_2(t)=\frac{1}{\sqrt{2H}}(e^{t\vec H})^*[2H\vec h_0-h_0\vec H],\\
    f_3(t)=-\frac{1}{\sqrt{2H}}\vec H,
  \end{array}
\right.
\]
where $\kappa=v_1c_{12}^2-v_2c_{12}^1-(c_{12}^1)^2-(c_{12}^2)^2-c_{01}^2$.
\end{thm}

If we assume that the flow of the Reeb field $v_0$ defines a free and proper group action, then the quotient $N$ of the manifold $M$ by this group action is a manifold and the subriemannian metric on $M$ induces a Riemannian metric on $N$. In this case, $\kappa$ is simply the Gauss curvature of $N$ (see Section \ref{proofSas} for the proof of the following proposition).

\begin{prop}\label{Gauss}
Assume that the Reeb field $v_0$ defines a proper $G$-action ($G=S^1$ or $\Real$) on the subriemannian manifold $M$. If the quotient manifold $N=M/G$ is equipped with the Riemannian metric induced by the subriemannian one on $M$. Then the Gauss curvature of $N$ coincides with $\kappa$ defined in Theorem \ref{structuraldetailSas}.
\end{prop}

In particular, Proposition \ref{Gauss} shows that $\mathbb H^3$, $SU(2)$, and $SL(2)$ with standard subriemannian structures defined in Section \ref{contact} satisfies $\kappa=0$, $\kappa=1$, and $\kappa=-1$, respectively.

\smallskip

\section{Sasakian Space Forms and Generalized Measure Contraction Property}\label{GMCP}

In this section, we specialize the definition of measure contraction property to the three dimensional contact subriemannian case and rewrite it as a condition on the volume growth of the Popp's measure along subriemannian geodesics. Then we go on and compute explicitly this volume growth for the Sasakian manifolds with $\kappa$ defined in Theorem \ref{structuraldetail} equal to a constant. We will refer to these Sasakian manifolds as Sasakian space forms. With this as a motivation, we will introduce the generalized measure contraction property $\mathcal{MCP}(K;2,3)$ at the end.

Let $(M,\Delta,\left<\cdot,\cdot\right>)$ be a contact subriemannian manifold with subriemannian distance function $d$ and let $x_0$ be a point in $M$. Let $\mathfrak f$ be the function defined by $\mathfrak f(x)=-\frac{1}{2}d^2(x_0,x)$. According to the result in \cite{AgLe1}, the function $\mathfrak f$ is Lipschitz with respect to a Riemannian distance. In particular, it is differentiable almost everywhere. Therefore, we can define the map $\varphi_t$ by
\begin{equation}\label{alonggeo}
\varphi_t(x)=\pi(e^{t\vec H}(d\mathfrak f_x)),
\end{equation}
where $e^{t\vec H}$ is the subriemannian geodesic flow and $\pi:T^*M\to M$ is the natural projection.

For each fixed $x$ in the contact subriemannian manifold $M$, the curve $t\mapsto\varphi_t(x)$ is a minimizing geodesic starting from $x$ and ending at $x_0$. In particular, $\varphi_1$ is the constant map $\varphi_1(x)=x_0$. Moreover, since the function $\mathfrak f$ is Lipschitz with respect to a Riemannian distance, $t\mapsto\varphi_t(x)$ is uniquely minimizing between its end-points for Lebesgue almost all points $x$ (see \cite{AgLe1}). It follows that $\varphi_1$ is the unique solution to the optimal transportation problem (\ref{optimalorigin}) when the final measure $\mu_1$ is a delta mass $\delta_{x_0}$ at the point $x_0$. It also follows that the path of measures $\varphi_{t*}\mu$ defines a Wasserstein geodesic for any given measure $\mu$ which is absolutely continuous with respect to the Popp measure. Moreover, this is the only geodesic connecting $\mu$ and $\delta_{x_0}$. It follows from Definition \ref{MCP} and Remark \ref{switch} below that the measure contraction property is a control on the volume growth $\eta(\varphi_t(U))$ of the set $U$ along geodesics $t\mapsto\varphi_t(x)$ which end at $x_0$. In the case of Sasakian space forms, the volume growth $\eta(\varphi_t(U))$ is given by the following equality (see Section \ref{proofSasaVol} for the proof).

\begin{thm}\label{SasaVol}
Let $(M,\Delta,\left<\cdot,\cdot\right>)$ be a three dimensional Sasakian manifold with $\kappa=K$ a constant. Let $d$ be the subriemannian distance and $\eta$ be the Popp's measure. Let $x_0$ be a point on the manifold $M$ and let $\varphi_t$ be defined as in (\ref{alonggeo}). Then the following holds
\[
\eta(\varphi_t(U))=\int_U(1-t) \left(\frac{s(\mathfrak k(x),(1-t)D(x))}{s(\mathfrak k(x),D(x))}\right) d\eta(x)
\]
for any Borel set $U$, where
\[
s(k,r)=
\begin{cases}
\frac{12(2-2\cos(\sqrt kr)-\sqrt k\,r\sin(\sqrt k\, r)}{k^2} & \text{if } k>0\\
r^4 & \text{if } k=0\\
\frac{12(2-2\cosh(\sqrt{-k}\,r) +(\sqrt{-k}\,r)\sinh(\sqrt{-k}\,r))}{k^2} & \text{if } k<0
\end{cases}
\]
$\mathfrak k(x)=(v_0D)^2(x)+K$, and $D(x)=d(x_0,x)$.
\end{thm}

Note that
\[
\frac{s(\mathfrak k(x),(1-t)D(x))}{s(\mathfrak k(x),D(x))}\geq \frac{s(K,(1-t)D(x))}{s(K,D(x))}.
\]

In view of this and Theorem \ref{SasaVol}, we define the generalized measure contraction property as follows.

\begin{defn}\label{genMCP}
The metric measure space $(M,d,\eta)$ satisfies the \textit{generalized measure contraction property} $\mathcal{MCP}(K;2,3)$ if for each $t$ in $(0,1)$, there is a Markov kernel $P_t$ which takes a point in $\mathfrak X\times\mathfrak X$ to a measure in $\mathfrak X$ such that, for $\nu^2$ almost every $(x_0, x)$ and for $P_t(x_0, x)$ almost every $z$, the point $z$ is a $t$-intermediate point of $x_0$ and $x$, and the followings hold:
\begin{equation}\label{gmcp1}
\eta(U)\geq  \int_M(1-t) \left(\frac{\mathfrak s_K((1-t)D(x))}{\mathfrak s_K(D(x))}\right)P_t(x,x_0)(U) d\eta(x)
\end{equation}
and
\begin{equation}\label{gmcp2}
\eta(U)\geq  \int_Mt \left(\frac{\mathfrak s_K(tD(x))}{\mathfrak s_K(D(x))}\right)P_t(x_0,x)(U) d\eta(x)
\end{equation}
for any measurable set $U$ and $\nu$-almost every $x$ in $M$, where
\[
\mathfrak s_K(r)=
\begin{cases}
\frac{12(2-2\cos(\sqrt Kr)-\sqrt K\,r\sin(\sqrt K\, r)}{K^2} & \text{if }K>0\\
r^4 & \text{if }K=0\\
\frac{12(2-2\cosh(\sqrt{-K}\,r) +(\sqrt{-K}\,r)\sinh(\sqrt{-K}\,r))}{K^2} & \text{if }K<0
\end{cases}
\]
and $D(x)=d(x_0,x)$.
\end{defn}

\begin{rem}\label{switch}
In the subriemannian case, if $t\mapsto\gamma(t)$ is a minimizing geodesic satisfying $\gamma(0)=x$ and $\gamma(1)=y$, then $t\mapsto\gamma(1-t)$ is a minimizing geodesic going from $y$ to $x$. It follows from this that (\ref{gmcp1}) implies (\ref{gmcp2}) in the subriemannian case. Moreover, if $\eta$ is absolutely continuous with respect to the Popp volume and the subriemannian manifold is contact, then $P_t(x,x_0)=\delta_{\varphi_t(x)}$ for $\eta$ almost every $x$. Therefore, (\ref{gmcp1}) becomes
\[
\eta(U)\geq  \int_{\varphi_t^{-1}(U)}(1-t) \left(\frac{\mathfrak s_K((1-t)D(x))}{\mathfrak s_K(D(x))}\right) d\eta(x).
\]
On the other hand, if the following holds instead
\begin{equation}\label{gmcp3}
\eta(\varphi_t(B))\geq  \int_{B}(1-t) \left(\frac{\mathfrak s_K((1-t)D(x))}{\mathfrak s_K(D(x))}\right) d\eta(x)
\end{equation}
for any measurable set $B$, then
\[
\begin{split}
&\eta(U)\geq \eta(\varphi_t(\varphi_t^{-1}(U)))\\
&\geq\int_{\varphi_t^{-1}(\varphi_t(\varphi_t^{-1}(U)))}(1-t) \left(\frac{\mathfrak s_K((1-t)D(x))}{\mathfrak s_K(D(x))}\right) d\eta(x)\\
&\geq\int_{\varphi_t^{-1}(U)}(1-t) \left(\frac{\mathfrak s_K((1-t)D(x))}{\mathfrak s_K(D(x))}\right) d\eta(x).
\end{split}
\]
Therefore, it is enough to verify (\ref{gmcp3}) in order to verify $\mathcal{MCP}(K;2,3)$.
\end{rem}

\begin{rem}
If $\Delta$ is a bracket-generating distribution, then it defines a flag of distribution by
\[
\Delta^1:=\Delta\subset\Delta^2\subset...\subset TM.
\]
If we denote the dimension of the vector space $\Delta^i_x$ by  $n^i_x$, then
the growth vector of the distribution $\Delta$ at the point $x$ is defined by
\[
(n^1_x,n^2_x,...,n^k_x).
\]
The pair (2, 3) in the generalized measure contraction
property is the growth vector of the three dimensional contact subriemannian
manifold. In this paper, we add $\mathcal{MCP}(K; 2,3)$ to the measure
contraction property $MCP(K,N)$ introduced earlier by Sturm. It would
be very interesting to find appropriate measure contraction properties
for other subriemannian manifolds with different growth vectors.
\end{rem}

\begin{rem}
Note that the condition $\mathcal{MCP}(0;2,3)$ is the same as $MCP(0,5)$.
\end{rem}

\begin{rem}
$\mathfrak s_K$ in the Definition \ref{genMCP} satisfies
\[
\mathfrak s_K(r)=r^4+o(r^4) \text{ as } r\to 0.
\]
Therefore, $\mathcal{MCP}(K;2,3)$ does not imply $MCP(0,N)$ for any $N<5$.
\end{rem}

\smallskip

\section{The Main Result and its Consequences}\label{main}

In this section, we state our main result and its consequences. For their proofs, see Section \ref{proofmain}.

\begin{thm}\label{mainSasa}(Generalized Measure Contraction Property)
Assume that the three dimensional contact subriemannian manifold $M$ is Sasakian (i.e. $R^{22}\equiv 0$). Then the followings are equivalent:
\begin{itemize}
\item there is a constant $K$ such that $\ric(\alpha)\geq 2KH(\alpha)$  for all $\alpha$ in the cotangent bundle $T^*M$,
\item $\kappa$ is bounded below by $K$,
\item the metric measure space $(M,d,\eta)$ satisfies the generalized measure contraction property $\mathcal{MCP}(K;2,3)$, where $d$ is the subriemannian distance and $\eta$ is the Popp's measure (see Section \ref{subgeo} for the definitions).
\end{itemize}
\end{thm}

Recall that $\mathcal{MCP}(0;2,3)$ is the same as $MCP(0,5)$. Therefore, Theorem \ref{SasaMCP} follows from Theorem \ref{mainSasa}.

\begin{rem}
As mentioned in Remark \ref{TW}, $\kappa$ is the Tanaka-Webster curvature. Therefore, Theorem \ref{mainSasa} provides an alternative characterization of when the Tanaka-Webster curvature of a Sasakian manifold is bounded below by a constant $K$.
\end{rem}

\begin{rem}\label{Igor}
The proof of Theorem \ref{structuraldetail} also works when we assume that $R^{22}_\alpha\geq 0$ for all $\alpha$ in the cotangent bundle. However, it is I. Zelenko's observation (private communications) that $R^{22}_\alpha\geq 0$ for all $\alpha$ implies $R^{22}\equiv 0$. On the other hand, see Section \ref{more} for result with relaxed assumption on $R^{22}$.
\end{rem}

\begin{rem}
Many ingredients used in the proof of Theorem \ref{mainSasa} are also present in the higher dimensional contact subriemannian case. This includes the recent result in \cite{LiZe1,LiZe2}, a comparison principle of matrix Riccati equations, and the solvability of matrix Riccati equations with constant coefficients. Therefore, results similar to Theorem \ref{mainSasa} can be proved in a similar way in the higher dimensional case where the canonical Darboux frames and curvature invariants are well understood (i.e. an analog of Theorem \ref{structuraldetail}). For instance, the result in \cite{Ju} for the higher dimensional Heisenberg group can be proved in the same way as in Theorem \ref{mainSasa}.
\end{rem}

Let $B_x(R)$ be the subriemannian ball of radius $R$ centered at a point $x$ in the manifold $M$ and let $\pi:T^*M\to M$ be the natural projection. The proof of Theorem \ref{mainSasa} is still valid if the curvature assumptions only holds on a ball $B_x(R)$ and the measure is contracted towards the center of the ball $x$. Therefore, the following volume doubling property holds.

\begin{cor}\label{doublingSasa}(Volume Doubling Property)
Assume that there is a point $x_0$ in the three dimensional contact subriemannian manifold and a constant $R>0$ such that  $R^{11}_\alpha\geq 0$ and $R^{22}_\alpha=0$ (i.e. $M$ is Sasakian) for all $\alpha$ in $\pi^{-1}(B_{x_0}(2R))$ and  for some constant $K\geq 0$. Then
\[
 \eta(B_{x_0}(2kR))\leq 2^5\eta(B_{x_0}(kR))
\]
for all $0<k<1$.
\end{cor}

\begin{rem}
Note that although the generalized measure contraction property is sharp (see Section \ref{GMCP}), the constant $2^5$ in Corollary \ref{doublingSasa} is not. This is very different from the Riemannian case and it is due to a key difference between the Riemannian and subriemannian cut locus. Given a point $x$ in a Riemannian manifold, there is a small enough neighborhood containing $x$ which does not contain any cut point of $x$. On the other hand, any neighborhood of a point $x$ has a nonempty intersection with the cut locus of $x$ in the subriemannian case (see \cite{Ag2}). In particular, we don't obtain a family of shrinking balls if we contract a subriemannian ball along geodesics to the center of the ball. This is very different from the Riemannian case. For the sharp constant in Corollary \ref{doublingSasa}, see \cite{AgLe2}.
\end{rem}

The local Poincar\'e inequality also holds under the assumptions in Corollary \ref{doublingSasa}. For this, let $\nabla_Hf$ be the \textit{horizontal gradient} of the function $f$ defined by the condition $df(v)=\left<\nabla_H f,v\right>$ for all $v$ in the distribution $\Delta$. For the proof of the following corollary, see Section \ref{proofmain}.

\begin{cor}\label{PoincareSasa}(Local Poincar\'e Inequality)
Under the assumptions in Theorem \ref{mainSasa}, the following local Poincar\'e inequality holds for all smooth functions $f$ and all $0<k<1$
\[
\begin{split}
&\frac{1}{\eta(B_{x_0}(kR))}\int_{B_{x_0}(kR)}|f(x)-\left<f\right>_{B_{x_0}(kR)}|d\eta(x) \\
&\leq \frac{CR}{\eta(B_{x_0}(2kR))}\int_{B_{x_0}(2kR)}|\nabla_Hf|d\eta(x),
\end{split}
\]
for some constant $C$ and where
\[
\left<f\right>_{B_{x_0}(kR)}=\frac{1}{\eta({B_{x_0}(kR)})}\int_{B_{x_0}(kR)}f(x)d\eta(x).
\]
\end{cor}

Let $\Delta_H$ be the sub-Laplacian defined by $\Delta_H=\textbf{div}_\eta\nabla_H$, where $\textbf{div}_\eta$ denotes the divergence with respect to $\eta$. Under the assumptions in Theorem \ref{SasaMCP}, the results in \cite{CoHoSa} together with Corollary \ref{doublingSasa} and \ref{PoincareSasa} show that any positive harmonic function of the sub-Laplacian $\Delta_H$  satisfies the Harnack inequality. More precisely,

\begin{thm}(Harnack inequality for sub-Laplacian)\label{Harnack}
Under the assumptions in Corollary \ref{doublingSasa}, any positive solution to the equation $\Delta_Hf=0$ satisfies
\[
\sup_{B_{x_0}(kR)}f\leq C\inf_{B_{x_0}(kR)}f
\]
for all $0<k<1$.
\end{thm}

For the proof of Theorem \ref{Harnack}, see \cite{CoHoSa}. Finally, by letting $R$ goes to $+\infty$ in Theorem \ref{Harnack}, the following Liouville theorem holds.

\begin{cor}(Liouville Theorem for sub-Laplacian)
Under the assumptions in Theorem \ref{SasaMCP}, any non-negative solution to the equation $\Delta_Hf=0$ is a constant.
\end{cor}

In the special case when the manifold $M$ is compact, the above Harnack inequality and Liouville Theorem were done in \cite{CaYa}.

\smallskip

\section{More General Situations and Final Remark}\label{more}

In this section, we show that the assumption on $R^{22}$ in Theorem \ref{mainSasa} can be relaxed. To do this, let $\Omega_x$ be the injectivity domain at a point $x$ in $M$ defined as the set of all covectors $\alpha$ in $T^*_xM$ such that
\[
t\mapsto \pi(e^{t\vec H}(\alpha)), 0\leq t\leq 1
\]
is length minimizing between its end points. Finally, let $\Omega=\bigcup_x\Omega_x$ be the injectivity domain.

One can apply similar arguments as in the proof of Theorem \ref{mainSasa} under the assumption that $R^{22}$ is bounded below by a constant on $\Omega$ instead of bounded by zero. This will give certain measure contraction property.

\begin{thm}\label{morethm}
Assume that $M$ is a three dimensional contact subriemannian manifold with subriemannian distance $d$ and Popp's measure $\eta$. Assume further that there is a constant $C_1$ and a non-negative constant $C_2$ such that $R^{11}_\alpha\geq 2C_1H(\alpha)$ and $R^{22}_\alpha\geq -C_2^2$ for all $\alpha$ in $\Omega$. Let $\varphi_t$ be as in (\ref{alonggeo}). Then the metric measure space $(M,d,\eta)$ satisfies
\[
\eta(\varphi_t(U))\geq  \int_U(1-t) \left(\frac{s_{C_1,C_2}(\sqrt 2 (1-t))}{s_{C_1,C_2}(\sqrt 2)}\right) d\eta(x)
\]
for any Borel set $U$, where
\[
s_{C_1,C_2}(r)=\begin{cases}
\frac{\cosh(r\sqrt\mathfrak a)-1}{\mathfrak a}+\frac{\cos(r\sqrt\mathfrak b)-1}{\mathfrak b} & \text{if } \mathfrak a>0 \text{ and }\mathfrak b>0,\\
\frac{r^2}{2}+\frac{\cos(r\sqrt\mathfrak b)-1}{\mathfrak b} & \text{if } \mathfrak a=0 \text{ and }\mathfrak b>0,\\
\frac{\cosh(r\sqrt\mathfrak a)-1}{\mathfrak a}-\frac{r^2}{2} & \text{if } \mathfrak a>0 \text{ and }\mathfrak b=0,\\
\frac{\cosh(r\sqrt\mathfrak a)-1}{\mathfrak a}+\frac{\cosh(r\sqrt{-\mathfrak b})-1}{\mathfrak b} & \text{if } \mathfrak a>0 \text{ and }\mathfrak b<0,\\
\frac{r^2}{2}+\frac{\cosh(r\sqrt{-\mathfrak b})-1}{\mathfrak b} & \text{if } \mathfrak a=0 \text{ and }\mathfrak b<0,\\
\frac{\cos(r\sqrt{-\mathfrak a})-1}{\mathfrak a}+\frac{\cos(r\sqrt\mathfrak b)-1}{\mathfrak b} & \text{if } \mathfrak a<0 \text{ and } \mathfrak b>0,\\
\frac{\cos(r\sqrt\mathfrak a)-1}{\mathfrak a}-\frac{r^2}{2} & \text{if } \mathfrak a<0 \text{ and }\mathfrak b=0,\\
r^4 & \text{if } \mathfrak a=\mathfrak b=0,
\end{cases}
\]
$\mathfrak a(x)=C_2-\frac{1}{2}C_1d^2(x_0,x)$, and $\mathfrak b(x)=C_2+\frac{1}{2}C_1d^2(x_0,x)$.
\end{thm}

The proof of Theorem \ref{morethm} is very similar to that of Theorem \ref{mainSasa}. We only outline the differences here and omit the detailed proof of Theorem \ref{morethm}. In the proof of Theorem \ref{mainSasa}, we use a comparison theorem of matrix Riccati equations to compare (\ref{Riccati}) with (\ref{sriccati}). Since $R^{22}_\alpha$ in Theorem \ref{morethm} is bounded below instead of vanishes, we have to change $\tilde R_\alpha$ in (\ref{sriccati}) to
\[
\tilde R_\alpha=\left(
\begin{array}{ccc}
           2C_1H(\alpha) & 0 & 0 \\
           0 & -C_2^2 & 0 \\
           0 & 0 & 0 \\
         \end{array}\right)
\]
for the proof of Theorem \ref{morethm}. The resulting equation (\ref{sriccati}) is still a Riccati equation with constant coefficient and can therefore be integrated. The rest of the proof of Theorem \ref{morethm} is the same as that of Theorem \ref{mainSasa}

Finally, we show that any compact three dimensional contact subriemannian manifold satisfies the assumptions in Theorem \ref{morethm}.

\begin{thm}\label{injcompact}
Assume that the three dimensional contact subriemannian manifold is compact. Then $R^{22}_\alpha\Big|_{\alpha\in \Omega}$ is bounded. In particular, it satisfies the assumptions in Theorem \ref{morethm}.
\end{thm}

For the proof of Theorem \ref{injcompact}, see Section \ref{proofinj}.

\smallskip

\section{Proof of Theorem \ref{structural}, \ref{structuraldetail}, and \ref{R22}}\label{proofstructural}

In this section, we give the proof of Theorem \ref{structural}, \ref{structuraldetail}, and \ref{R22}. Let us start with a lemma on Euler field. Recall that $\vec E$ denotes the Euler field and $H$ denotes the subriemannian Hamiltonian.

\begin{lem}\label{homogeneous}
$(e^{t\vec H})^*\vec E=\vec E-t\vec H$
\end{lem}

\begin{proof}
Recall $\mathfrak d_s:T^* M\to T^* M$ is the dilation map $\mathfrak d_s(\alpha)=s\alpha$. By the definition of the symplectic form,
\[
\mathfrak d_s^*\omega=s\omega.
\]

It follows that
\[
\begin{split}
&\omega(d\mathfrak d_s(\vec H(\alpha)),X(s\alpha))\\
&=s\omega(\vec H(\alpha),d\mathfrak d_s^{-1}(X(s\alpha)))\\
&=-sdH(d\mathfrak d_{1/s}(X(s\alpha))),
\end{split}
\]
where $X$ is any tangent vector in the tangent bundle $TT^*M$.

The subriemannian Hamiltonian $H$ is homogeneous of degree two in the fibre direction. In other words,
\[
H(\mathfrak d_s(\alpha))=s^2H(\alpha).
\]
Therefore,
\[
\omega(d\mathfrak d_s(\vec H(\alpha)),X(s\alpha))=-\frac{1}{s}dH(X(s\alpha))=\frac{1}{s}\omega(\vec H(s\alpha),X(s\alpha)).
\]

It follows that $\mathfrak d_s^*\vec H=s\vec H$, where $\mathfrak d_s^*\vec H$ is the pullback of the vector field $\vec H$ by the map $\mathfrak d_s$. By comparing the flow of the above vector fields, we have
\[
e^{t\vec H}\circ\mathfrak d_s=\mathfrak d_s\circ e^{ts\vec H}.
\]
By differentiating the above equation with respect to $s$ and set $s$ to 1, it follows that $(e^{t\vec H})^*\vec E=\vec E-t\vec H$ as claimed.
\end{proof}

\begin{proof}[Proof of Theorem \ref{structural}]
According to the main result in \cite{LiZe1,LiZe2}, there exists a family of Darboux frames
\[
\{e_1(t),e_2(t),e_3(t),f_1(t),f_2(t),f_3(t)\}
\]
and functions $R^{ij}_\alpha(t)$ which satisfy
\[
\left\{
  \begin{array}{ll}
    \dot e_1(t)=f_1(t), \\
    \dot e_2(t)=e_1(t), \\
    \dot e_3(t)=f_3(t), \\
    \dot f_1(t)=-R^{11}_\alpha(t)e_1(t)-R^{31}_\alpha(t)e_3(t)-f_2(t),\\
    \dot f_2(t)=-R^{22}_\alpha(t)e_2(t)-R^{32}_\alpha(t)e_3(t),\\
    \dot f_3(t)=-R^{31}_\alpha(t)e_1(t)-R^{32}_\alpha(t)e_2(t)-R^{33}_\alpha(t)e_3(t).
  \end{array}
\right.
\]

\begin{rem}
In the language of \cite{LiZe1,LiZe2}, the Young diagram associated with the above structural equations consists of two columns with two boxes in the first column and one box in the second column. Note that the reduced and the non-reduced Young diagrams are the same in this case.
\end{rem}

Note that $d\pi(\vec E)=0$. Therefore, $\vec E(e^{t\vec H}(\alpha))$ is contained in the vertical space at $e^{t\vec H}(\alpha)$ for each time $t$. Hence, by the definition of the Jacobi curve $J_\alpha(t)$, the vector $(e^{t\vec H})^*\vec E(\alpha)$ is contained in $J_\alpha(t)$ for each $t$. It follows from Lemma \ref{homogeneous} that
\[
\vec{E}(\alpha)-t\vec{H}(\alpha)=\sum\limits_{i=1}^3a_i(t)e_i(t)
\]
for some functions $a_i$ of time $t$. If we differentiate with respect to time $t$ twice, we get
\[
\begin{array}{ll}
2\dot a_1(t) f_1(t)+2\dot a_2(t)e_1(t)+2\dot a_3 (t)f_3(t)-a_1(t)(R^{11}_\alpha(t)e_1(t)+\\+R^{31}_\alpha(t)e_3(t)+f_2(t))+a_2(t) f_1(t)-a_3(t)(R^{31}_\alpha(t)e_1(t)+R^{32}_\alpha(t)e_2(t)+\\+R^{33}_\alpha(t)e_3(t))+\ddot a_1(t)e_1(t)+\ddot a_2(t)e_2(t)+\ddot a_3(t)e_3(t)=0.
\end{array}
\]
If we equate the coefficients of the $f_i(t)$'s, we get $a_1\equiv a_2\equiv \dot a_3\equiv 0$. Therefore, $\vec{E}(\alpha)-t\vec{H}(\alpha)=a_3e_3(t)$ and $-\vec{H}(\alpha)=a_3f_3(t)$ for some constant $a_3$ satisfying $(a_3)^2=\omega(a_3f_3(t),a_3e_3(t))=dH(\vec{E}(\alpha))=2H(\alpha)$. It follows that $R^{31}_\alpha(t)=R^{32}_\alpha(t)=R^{33}_\alpha(t)=0$. Moreover, we also have
\begin{equation}\label{e3f3}
e_3(t)=\frac{1}{(2H(\alpha))^{1/2}}(\vec{E}(\alpha)-t\vec{H}(\alpha)),\quad f_3(t)=-\frac{1}{(2H(\alpha))^{1/2}}\vec{H}(\alpha).
\end{equation}
\end{proof}

For the proof of Theorem \ref{structuraldetail}, we need a few more lemmas. Recall that $h_{ij}:T^*M\to\Real$ be the Hamiltonian lift of the vector field $[v_i,v_j]$ defined by
\[
h_{ij}(\alpha)=\alpha([v_i,v_j]).
\]
The commutator relations of the frame $\{\vec h_i,\vec\alpha_i|i=0,1,2\}$ are given by the following:

\begin{lem}\label{bundlemain1}
\[
\begin{split}
&[\vec h_i,\vec h_j]=\vec h_{ij},\quad [\vec h_i,\vec \alpha_j]=-\sum\limits_k c_{ik}^j\vec\alpha_k,\quad [\vec\alpha_i,\vec\alpha_j]=0,\\
&[\vec h_1,\vec\alpha_0]=\vec\alpha_2,\quad [\vec h_2,\vec\alpha_0]=-\vec\alpha_1
\end{split}
\]
\end{lem}

\begin{proof}
Since the Lie derivative $\mathcal L$ of the symplectic form $\omega$ along the Hamiltonian vector field $\vec h_i$ vanishes,
\begin{equation}\label{hihj}
i_{[\vec h_i,\vec h_j]}\omega=\mathcal L_{\vec h_i}i_{\vec h_j}\omega-i_{\vec h_j}\mathcal L_{\vec h_i}\omega=\mathcal L_{\vec h_i}i_{\vec h_j}\omega=-d(\omega(\vec h_i,\vec h_j)).
\end{equation}

The function $\omega(\vec h_i,\vec h_j)$ is equal to $h_{ij}$. Indeed, since $d\pi(\vec h_i)=v_i$, we have
\[
\theta_\alpha(\vec h_i)=\alpha(d\pi(\vec h_i))=\alpha(v_i)=h_i(\alpha).
\]
It follows from this and the Cartan's formula that
\[
\begin{split}
&dh_j(\vec h_i)=\omega(\vec h_i,\vec h_j)=d\theta(\vec h_i,\vec h_j)\\
&=\vec h_i(\theta(\vec h_j))-\vec h_j(\theta(\vec h_i))-\theta([\vec h_i,\vec h_j])\\
&=dh_j(\vec h_i)-dh_i(\vec h_j)-\theta([\vec h_i,\vec h_j]).
\end{split}
\]
If we apply again $d\pi(\vec h_i)=v_i$, then we have
\[
\theta_\alpha([\vec h_i,\vec h_j])=\alpha(d\pi([\vec h_i,\vec h_j]))=\alpha([v_i,v_j])=h_{ij}(\alpha).
\]
Therefore, we have
\begin{equation}\label{hihj1}
\omega(\vec h_i,\vec h_j)=-dh_i(\vec h_j)=h_{ij}.
\end{equation}
If we combine this with (\ref{hihj}), the first assertion of the lemma follows.

A calculation similar to the above one shows that
\[
i_{[\vec h_i,\vec \alpha_j]}\omega=\mathcal L_{\vec h_i}i_{\vec \alpha_j}\omega.
\]
By Cartan's formula, the above equation becomes
\[
i_{[\vec h_i,\vec \alpha_j]}\omega=-i_{\vec h_i}\pi^*d\alpha_j=-\pi^*(i_{v_i}d\alpha_j).
\]
The second assertion follows from this and (\ref{dualbracket}).

If we apply Cartan's formula again,
\[
i_{[\vec\alpha_i,\vec\alpha_j]}\omega
=\mathcal L_{\vec\alpha_i}i_{\vec\alpha_j}\omega-i_{\vec\alpha_j}\mathcal L_{\vec\alpha_i}\omega
=-i_{\vec\alpha_i}d(\pi^*\alpha_j)+i_{\vec\alpha_j}d(\pi^*\alpha_i)
\]
Since $d\pi(\vec\alpha_i)=0$, it follows that $i_{[\vec\alpha_i,\vec\alpha_j]}\omega=0$. Therefore, the third assertion holds by the non-degeneracy of $\omega$.

Finally, the last two assertions follows from Lemma \ref{strcst}.
\end{proof}

Let $\beta=h_1dh_2-h_2dh_1$, then we also have the following relations:

\begin{lem}\label{bundlemain2}
\[
dh_i(\vec h_j)=-h_{ij},\quad  \alpha_i(\vec h_j)=-dh_i(\vec\alpha_j)=\delta_{ij},\quad \alpha_i(\vec\alpha_j)=0,
\]
\[
 \beta(\vec\xi_2)=dH(\vec\xi_1)=0,\quad \beta(\vec\xi_1)=dH(\vec\xi_2)=-2H, \quad \beta(\vec H)=2Hh_{12}
\]
\end{lem}

\begin{proof}
The first assertion follows from (\ref{hihj1}) and the next two assertions follow from $d\pi(\vec h_i)=v_i$ and $d\pi(\vec\alpha_i)=0$. A computation using $\alpha_i(\vec h_j)=\delta_{ij}$ proves the third and the fourth assertions. The final assertion follows from the following computations
\[
\beta(\vec H)=(h_1dh_2-h_2dh_1)(h_1\vec h_1+h_2\vec h_2)=h_1^2dh_2(\vec h_1)-h_2^2dh_1(\vec h_2)=2Hh_{12}.
\]
\end{proof}

\begin{proof}[Proof of Theorem \ref{structuraldetail}]
Recall $J_\alpha(\cdot)$ denotes the Jacobi curve at the point $\alpha$ in the cotangent bundle $T^*M$. By Theorem \ref{structural}, there exists a family of Darboux frame
\[
\{e_1(t),e_2(t),e_3(t),f_1(t),f_2(t),f_3(t)\}
\]
and functions $R^{ij}_\alpha(t)$ such that
\[
J_\alpha(t)=\textbf{span}\{e_1(t),e_2(t),e_3(t)\}
\]
and
\[
\left\{
  \begin{array}{ll}
    \dot e_1(t)=f_1(t), \\
    \dot e_2(t)=e_1(t), \\
    \dot e_3(t)=f_3(t), \\
    \dot f_1(t)=-R^{11}_\alpha(t)e_1(t)-f_2(t),\\
    \dot f_2(t)=-R^{22}_\alpha(t)e_2(t),\\
    \dot f_3(t)=0.
  \end{array}
\right.
\]

Let $\mathcal E(t)$ be defined by
\[
\mathcal E(t)=(e^{t\vec H})^*\vec\alpha_0(\alpha)=de^{-t\vec H}(\vec\alpha_0(e^{t\vec H}(\alpha))).
\]
By the definition of the Jacobi curve $J_\alpha(\cdot)$, we known that $\mathcal E(t)$ is contained in $J_\alpha(t)$ for each $t$. Since $e_1(t),e_2(t),e_3(t)$ span $J_\alpha(t)$, we must have
\[
\mathcal E(t)=c_1(t)e_1(t)+c_2(t)e_2(t)+c_3(t)e_3(t)
\]
for some functions $c_i$ of time $t$, $i=1,2,3$.

Let $\pi:T^*M\to M$ be the natural projection. The Hamiltonian vector field $\vec H$ of the subriemannian Hamiltonian $H$ satisfies
\[
d\pi(\vec H(\alpha))=h_1(\alpha)v_1+h_2(\alpha)v_2.
\]
It follows that
\[
\omega(\vec\alpha_0,\vec H)=-\pi^*\alpha_0(\vec H)=0.
\]

Since the flow $e^{t\vec H}$ preserves the symplectic form $\omega$, it follow from the definition of $\mathcal E(t)$ that
\[
 \omega(\mathcal E,\vec H)=0.
\]

By (\ref{e3f3}), we know that $f_3(t)=-\frac{1}{(2H)^{1/2}}\vec H$. Since $\{e_i(t),f_i(t)|i=1,2,3\}$ is a Darboux basis, we have
\[
 0=\omega(\mathcal E,\vec H)=(2H)^{1/2}c_3(t).
\]

This shows that $c_3\equiv 0$ and so
\[
\mathcal E(t)=c_1(t)e_1(t)+c_2(t)e_2(t).
\]
By the definition of $\mathcal E(t)$, if we differentiate this with respect to time $t$, then we have
\[
(e^{t\vec H})^*[\vec H,\vec\alpha_0]=\dot{\mathcal E}(t)=\dot c_1(t)e_1(t)+c_1(t)f_1(t)+\dot c_2(t)e_2(t)+c_2(t)e_1(t).
\]

By the Cartan's formula and $\alpha_0(\vec H)=0$, it follows that
\[
\omega(\dot{\mathcal E}(t),\mathcal E(t))=\omega([\vec H,\vec\alpha_0],\vec\alpha_0)=\pi^*\alpha_0([\vec H,\vec\alpha_0])=-\pi^*d\alpha_0(\vec H,\vec\alpha_0)=0.
\]

By combining this with the above equation for $\mathcal E$ and $\dot{\mathcal E}$, we have $c_1\equiv 0$. If we differentiate the equation
\[
\mathcal E(t)=c_2(t)e_2(t)
\]
with respect to time $t$ again, we get
\[
\begin{split}
&(e^{t\vec H})^*(ad_{\vec H}(\vec\alpha_0))=\dot c_2(t)e_2(t)+c_2(t)e_1(t)\\
&(e^{t\vec H})^*(ad_{\vec H}^2(\vec\alpha_0))=\ddot c_2(t)e_2(t)+2\dot c_2(t)e_1(t)+c_2(t)f_1(t).
\end{split}
\]
Here $ad_{\vec H}$ denotes $ad_{\vec H}(\cdot)=[\vec H,\cdot]$.

Since $\{e_i(t),f_i(t)|i=1,2,3\}$ is a Darboux basis and the flow $e^{t\vec H}$ preserves the symplectic form $\omega$,
\[
\begin{split}
(c_2(t))^2&=\omega_\alpha((e^{t\vec H})^*ad_{\vec H}^2(\vec\alpha_0),(e^{t\vec H})^*ad_{\vec H}(\vec\alpha_0)))\\
&=(e^{t\vec H}_*\omega)_{e^{t\vec H}(\alpha)}(ad_{\vec H}^2(\vec\alpha_0),ad_{\vec H}(\vec\alpha_0)))\\
&=\omega_{e^{t\vec H}(\alpha)}(ad_{\vec H}^2(\vec\alpha_0),ad_{\vec H}(\vec\alpha_0)))\\
\end{split}
\]
Therefore, $c_2(t)=(e^{t\vec H})^*\left(\frac{1}{c}\right)$, where $c(\alpha):=\frac{1}{(\omega_\alpha(ad_{\vec H}^2(\vec\alpha_0),ad_{\vec H}(\vec\alpha_0)))^{1/2}}$.

It follows from the definition of $\mathcal E$ that
\[
e_2(t)=\frac{1}{c_2(t)}\mathcal E(t)=(e^{t\vec H})^*(c\vec\alpha_0).
\]

To find out what $c$ is more explicitly, we first compute $[\vec H,\vec\alpha_0]$. The Lie bracket is a derivation in each of its entries, so

\[
\begin{split}
[\vec H,\vec\alpha_0]&=[h_1\vec h_1+h_2\vec h_2,\vec\alpha_0]\\
&=-dh_1(\vec\alpha_0)\vec h_1-dh_2(\vec\alpha_0)\vec h_2+h_1[\vec h_1,\vec\alpha_0]+h_2[\vec h_2,\vec\alpha_0].
\end{split}
\]

It follows from this, Lemma \ref{bundlemain1}, and Lemma \ref{bundlemain2} that

\[
[\vec H,\vec\alpha_0]=h_1\vec\alpha_2-h_2\vec\alpha_1=\vec\xi_1.
\]

Next, we want to compute $[\vec H,\vec\xi_1]$. For this, let
\begin{equation}\label{hxiassume1}
[\vec H,\vec\xi_1]=k_0\vec\alpha_0+k_1\vec\xi_1+k_2\vec\xi_2+\sum_{i=0}^2\tilde c_i\vec h_i
\end{equation}
for some functions $\tilde c_i$ and $k_i$.

To compute $\tilde c_0$ for instance, we apply $\alpha_0$ on both sides of (\ref{hxiassume1}). Using Lemma \ref{bundlemain2} and Cartan's formula, we have $\tilde c_0=0$. Similar computation gives $\tilde c_1=-h_2$ and $\tilde c_2=h_1$. This shows that
\begin{equation}\label{hxiassume2}
[\vec H,\vec\xi_1]=k_0\vec\alpha_0+k_1\vec\xi_1+k_2\vec\xi_2+h_1\vec h_2-h_2\vec h_1.
\end{equation}

By applying $dh_0$ on both sides of (\ref{hxiassume2}) and using Lemma \ref{bundlemain2} again, we have $k_0=h_2h_{01}-h_1h_{02}+{\vec\xi_1}a$, where $a=dh_0(\vec H)$. Similar calculations using $\beta$ and $dH$ give

\begin{equation}\label{hxi}
[\vec H,\vec\xi_1]= h_1\vec h_2-h_2\vec h_1+\chi_0\vec\alpha_0+({\vec\xi_1} h_{12})\vec\xi_1-h_{12}\vec\xi_2.
\end{equation}
where $\chi_0=h_2h_{01}-h_1h_{02}+{\vec\xi_1}a$.

It follows that
\[
c^{-2}=\omega(ad^2_{\vec H}(\vec\alpha_0),ad_{\vec H}(\vec\alpha_0))=2H
\]
and $e_2(0)=\frac{1}{\sqrt{2H}}\vec\alpha_0$. It also follows from Theorem \ref{structural} that

\begin{equation}\label{ef}
\begin{array}{ll}
e_1(0)=\frac{1}{\sqrt{2H}}\vec\xi_1, \\
f_1(0)=\frac{1}{\sqrt{2H}}[\vec H,\vec\xi_1],\\
\dot f_1(0)=\frac{1}{\sqrt{2H}}[\vec H,[\vec H,\vec\xi_1]],\\
\ddot f_1(0)=\frac{1}{\sqrt{2H}}[\vec H,[\vec H,[\vec H,\vec\xi_1]]].
\end{array}
\end{equation}

A computation similar to that of (\ref{hxi}) gives

\begin{equation}\label{hhxi}
[\vec H,[\vec H,\vec\xi_1]]= -2H\vec h_0+h_0\vec H+\chi_1\vec\alpha_0+(\chi_2+\chi_0-{\vec\xi_1} a)\vec\xi_1+a\vec\xi_2
\end{equation}
where $\chi_1=h_0a+2{\vec H}{\vec\xi_1} a-{\vec\xi_1}{\vec H} a$ and
$\chi_2=h_0h_{12}+2{\vec H}{\vec\xi_1} h_{12}-{\vec\xi_1}{\vec H} h_{12}$.

It follows from Theorem \ref{structural}, (\ref{hxi}), (\ref{ef}) and (\ref{hhxi}) that
\begin{equation}\label{R11}
\begin{split}
R^{11}_\alpha(0)&=\omega(\dot f_1(0),f_1(0))\\
&=-\chi_0-\chi_2.
\end{split}
\end{equation}

Note that, in (\ref{R11}), $\vec\xi_1a$ does not appear. This is because
\[
\omega(-2H\vec h_0,h_1\vec h_2-h_2\vec h_1+\chi_0\vec\alpha_0)=-2H\vec\xi_1a
\]
and
\[
\omega(-(\vec\xi_1 a)\vec\xi_1, h_1\vec h_2-h_2\vec h_1)=2H\vec\xi_1a.
\]

Since $\dot f_1(0)=-R^{11}_\alpha(0)e_1(0)-f_2(0)$, it follows from (\ref{ef}), (\ref{hhxi}), and (\ref{R11}) that
\[
 f_2(0)=\frac{1}{\sqrt{2H}}[2H\vec h_0-h_0\vec H-\chi_1\vec\alpha_0+({\vec\xi_1} a)\vec\xi_1-a\vec\xi_2].
\]

A long computation using the bracket relations (\ref{bracket}) gives
\[
 \chi_2=-(h_0)^2+2H[(c_{12}^1)^2+(c_{12}^2)^2-{v_1}c_{12}^2+{v_2}c_{12}^1]+{\vec\xi_1}a.
\]
and
\[
\chi_0=h_2h_{01}-h_1h_{02}+{\vec\xi_1} a.
\]
We recall here that $v_ic_{jk}^l$ is the directional derivative of $c_{jk}^l$ in the direction $v_i$.

Another computation shows that
\[
\begin{split}
&\vec\xi_1 a=-\vec\xi_1(h_1h_{01}+h_2h_{02})\\
&=\vec\xi_1(c_{01}^1h_1^2+c_{01}^2h_1h_2+c_{02}^1h_1h_2+c_{02}^2h_2^2)\\
&=2c_{01}^1h_1dh_1(\vec\xi_1)+c_{01}^2dh_1(\vec\xi_1)h_2+c_{01}^2h_1dh_2(\vec\xi_1)\\
&+c_{02}^1dh_1(\vec\xi_1)h_2+c_{02}^1h_1dh_2(\vec\xi_1)+2c_{02}^2h_2dh_2(\vec\xi_1)\\
&=2c_{01}^1h_1h_2+c_{01}^2h_2^2-c_{01}^2h_1^2+c_{02}^1h_2^2-c_{02}^1h_1^2-2c_{02}^2h_1h_2\\
&=2h_{01}h_2-2h_{02}h_1-2Hc_{01}^2+2Hc_{02}^1.
\end{split}
\]

It follows as claimed that
\[
\begin{split}
&R^{11}_\alpha(0)=-\chi_0-\chi_2\\
&=-h_2h_{01}+h_1h_{02}-2{\vec\xi_1} a+(h_0)^2-2H[(c_{12}^1)^2+(c_{12}^2)^2-{v_1}c_{12}^2+{v_2}c_{12}^1]\\
&=-h_2h_{01}+h_1h_{02}-2{\vec\xi_1} a+(h_0)^2+2H\kappa+H(c_{01}^2-c_{02}^1)\\
&=h_0^2+2H\kappa-\frac{3}{2}{\vec\xi_1} a.
\end{split}
\]

To prove the formula for $R^{22}$, we differentiate the equation
\[
\dot f_1(t)=-R^{11}_\alpha(t) e_1(t)-f_2(t)
\]
and combine it with the equation
\[
\dot f_2(t)=-R^{22}_\alpha(t)e_2(t).
\]

We have
\[
R^{22}_\alpha(0)e_2(0)=\ddot f_1(0)+{\vec H} R^{11}_\alpha(0)e_1(0)+R^{11}_\alpha(0) f_1(0).
\]

Therefore, by applying $dh_0$ on both sides and using $dh_0(e_1(0))=0$, we get
\[
 R^{22}_\alpha(0)=-\sqrt{2H}[dh_0(\ddot f_1(0))+R^{11}_\alpha(0) dh_0(f_1(0))].
\]

By using Cartan's formula and (\ref{ef}), it follows that
\[
\begin{split}
&\sqrt{2H}dh_0(f_1(0))=dh_0([\vec H,\vec\xi_1])=-{\vec\xi_1} a,\\
&\sqrt{2H}dh_0(\dot f_1(0))=dh_0([\vec H,[\vec H,\vec\xi_1]])={\vec\xi_1}{\vec H} a-2{\vec H}{\vec\xi_1} a,\\
&\sqrt{2H}dh_0(\ddot f_1(0))=3{\vec H} {\vec\xi_1} {\vec H} a-3\vec H^2{\vec\xi_1} a- {\vec\xi_1} {\vec H}^2 a.
\end{split}
\]
The formula for $R^{22}_\alpha(0)$ follows from this.

\end{proof}

Finally, we come to the proof of Theorem \ref{R22}. The proof involves lengthy computations of $R^{22}$. Therefore, only a sketch is given below.

\smallskip

\begin{proof}[Proof of Theorem \ref{R22}]
Clearly, if $a\equiv 0$, then $R^{22}\equiv 0$ by Theorem \ref{structuraldetail}. Conversely, assume that $R^{22}\equiv 0$. By using the expression of $R^{22}$ in Theorem \ref{structuraldetail} and Lemma \ref{bundlemain2}, we can rewrite $R^{22}$ as a homogeneous polynomial of degree 4 with three variables $h_0$, $h_1$, and $h_2$. A long computation shows that the coefficients of $h_0^2h_1^2$ and $h_0^2h_1h_2$ are $-3(c_{01}^2+c_{02}^1)$ and $12c_{01}^1$, respectively. Therefore, if $R^{22}\equiv 0$, then $c_{01}^2+c_{02}^1=0$ and $c_{01}^1=-c_{02}^2=0$. It follows that
\[
\begin{split}
a&=dh_0(\vec H)\\
&=h_1dh_0(\vec h_1)+h_2dh_0(\vec h_2)\\
&=c_{02}^2h_1^2-(c_{01}^2+c_{02}^1)h_1h_2-c_{02}^2h_2^2\\
&=0.
\end{split}
\]
\end{proof}

\smallskip

\section{Proof of Theorem \ref{structuraldetailSas} and Proposition \ref{Gauss}}\label{proofSas}

In this section, we will give the proof of Theorem \ref{structuraldetailSas} and Proposition \ref{Gauss}. The result of Theorem \ref{structuraldetailSas} follows from the following two lemmas.

\begin{lem}\label{structuralconstant}
Under the assumptions of Theorem \ref{structuraldetailSas}, the functions $c_{ij}^k$ in the bracket relation (\ref{bracket}) satisfies
\[
c_{01}^1=c_{02}^2=0\quad\text{ and }\quad c_{01}^2=-c_{02}^1.
\]
\end{lem}

\begin{proof}[Proof of Lemma \ref{structuralconstant}]
If the flow of the vector field $v_0$ is denoted by $e^{tv_0}$, then the invariance of the subriemannian metric under the group action implies that
\[
\left<(e^{tv_0})^*v_i,(e^{tv_0})^*v_j\right>=\delta_{ij}\quad i,j=1,2.
\]
By differentiating the above equations with respect to time $t$, it follows that
\[
\alpha_j([v_0,v_i])+\alpha_i([v_0,v_j])=0 \quad i,j=1,2.
\]
If we apply the bracket relations (\ref{bracket}) of the frame $v_0,v_1,v_2$, we have
\[
c_{0i}^j+c_{0j}^i=\alpha_j([v_0,v_i])+\alpha_i([v_0,v_j])=0 \quad i,j=1,2.
\]
\end{proof}

It follows that

\begin{lem}\label{h0}
Under the assumptions of Theorem \ref{structuraldetailSas}, the function $h_0$ is a constant of motion of the flow $e^{t\vec H}$. i.e. $a=dh_0(\vec H)=0$.
\end{lem}

\begin{proof}[Proof of Lemma \ref{h0}]
This follows from general result in Hamiltonian reduction. In this special case this can also be seen as follow. By Lemma \ref{bundlemain2}

\begin{equation}\label{h0-1}
dh_0(\vec H)=dh_0(h_1\vec h_1+h_2\vec h_2)=h_1h_{10}+h_2h_{20}.
\end{equation}

By Lemma \ref{structuralconstant} we also have
\[
h_{10}=-c_{01}^0h_0-c_{01}^1h_1-c_{01}^2h_2=-c_{01}^2h_2.
\]

Similarly $h_{20}=-c_{02}^1h_1$. The result follows from this, (\ref{h0-1}), and Lemma \ref{structuralconstant}.
\end{proof}

\begin{proof}[Proof of Proposition \ref{Gauss}]

Let $\pi_M:M\to N$ be the quotient map. Let $w_1$ and $w_2$ be a local orthonormal frame on the surface $N$. Since $\pi$ is a submersion, there are unique vector fields $\tilde w_1$ and $\tilde w_2$ in the distribution $\Delta$ such that $d\pi(\tilde w_i)=w_i$. If $\Phi_t$ is the flow of the Reeb field $v_0$, then $\pi(\Phi_t(x))=\pi(x)$ by the definition of the quotient map. Therefore, $d\pi(d\Phi_t(\tilde w_i))=d\pi(w_i)$. Since $d\Phi_t(\tilde w_i)$ is in $\Delta$, we have $(\Phi_t)_*\tilde w_i=w_i$. If we differentiate this equation and set $t$ to zero, then we have $[v_0,\tilde w_i]=0$.

Since $\tilde w_1$ and $\tilde w_2$ are orthonormal with respect to the subriemannian metric, we can set $v_i=\tilde w_i$. It follows that $c_{01}^2=0$ and $\kappa$ is simplified to
\begin{equation}\label{magnetickappa}
\kappa={v_1}c_{12}^2-{v_2}c_{12}^1-(c_{12}^1)^2-(c_{12}^2)^2.
\end{equation}

From (\ref{bracket}), we also have $[v_1,v_2]=c_{12}^0v_0+c_{12}^1v_1+c_{12}^2v_2$. If we apply $d\pi_M$ to the equation, then we get $[w_1,w_2]=c_{12}^1w_1+c_{12}^2w_2$.

Let us denote the covariant derivative on the Riemannian manifold $N$ by $\nabla$. It follows from Koszul formula (\cite[Theorem 3.11]{On}) that
\begin{equation}\label{covariant}
\begin{array}{ll}
\nabla_{w_1}w_1=-c_{12}^1w_2,\quad \nabla_{w_2}w_2=-c_{12}^2w_1, \\
\quad \nabla_{w_1}w_2=c_{12}^1w_1, \quad \nabla_{w_2}w_1=-c_{12}^2w_2.
\end{array}
\end{equation}

Since the covariant derivative $\nabla$ is tensorial in the bottom slot and is a derivation in the other slot, it follows from (\ref{covariant}) that
\[
\begin{array}{ll}
\nabla_{[w_1,w_2]}w_1&=\nabla_{c_{12}^1w_1+c_{12}^2w_2}w_1\\ & =c_{12}^1\nabla_{w_1}w_1+c_{12}^2\nabla_{w_2}w_1\\ & =-[(c_{12}^1)^2+(c_{12}^2)^2]w_2
\end{array}
\]
and
\[
\begin{array}{ll}
[\nabla_{w_1},\nabla_{w_2}]w_1 &=\nabla_{w_1}\nabla_{w_2}w_1-\nabla_{w_2}\nabla_{w_1}w_1\\ &=-\nabla_{w_1}(c_{12}^2w_2)+\nabla_{w_2}(c_{12}^1w_2)\\ &=-({w_1}c_{12}^2)w_2+({w_2}c_{12}^1)w_2-2c_{12}^1c_{12}^2w_1.
\end{array}
\]

Therefore, it follows from the above calculation that the Gauss curvature is given by
\[
\begin{array}{ll}
<\nabla_{[w_1,w_2]}w_1-[\nabla_{w_1},\nabla_{w_2}]w_1,w_2>={w_1}c_{12}^2-{w_2}c_{12}^1-(c_{12}^1)^2-(c_{12}^2)^2.
\end{array}
\]

By (\ref{magnetickappa}), this agrees with $\kappa$.
\end{proof}

\smallskip

\section{Proof of Theorem \ref{SasaVol}}\label{proofSasaVol}

\begin{proof}[Proof of Theorem \ref{SasaVol}]
From the main result in \cite{CaRi}, the function $\mathfrak f(x)=-\frac{1}{2}d^2(x,x_0)$ is locally semi-concave on $M-\{x_0\}$, so it is differentiable almost everywhere. Assume that $x'$ is a point where $\mathfrak f$ is differentiable. It follows that the map $t\mapsto \varphi_t(x'):=\pi(e^{t\vec H}(d\mathfrak f_{x'}))$ is the unique minimizing geodesic connecting $x'$ and $x_0$. An argument similar to the Riemannian case using inverse function theorem shows that the function $\mathfrak f$ is $C^\infty$ in a neighborhood of the curve $t\mapsto \varphi_t(x')$ (see, for instance, \cite{Ja}). Moreover, it follows from \cite[Theorem 1.2]{Ag} that there is no conjugate point along the curve $t\mapsto\varphi_t(x')$. Therefore, the map $(d\varphi_t)_{x'}$ is nonsingular for each $t<1$.

If we denote the differential of the map $x\mapsto d\mathfrak f_x$ by $dd\mathfrak f$, then $d\varphi_t=d\pi (de^{t\vec H}(dd\mathfrak f))$. Let $e_i(t)$ and $f_i(t)$ be the Darboux frame at $d\mathfrak f_{x'}$ defined as in Theorem \ref{structural} and let $\varsigma_i=d\pi(f_i(0))$. Then the vectors $\{dd\mathfrak f(\varsigma_1),dd\mathfrak f(\varsigma_2),dd\mathfrak f(\varsigma_3)\}$ span a linear subspace $W$ of $T_{d\mathfrak f_{x'}} T^*M$. Therefore $dd\mathfrak f(\varsigma_i)$ can be written as
\begin{equation}\label{AB}
dd\mathfrak f(\varsigma_i)=\sum\limits_{j=1}^3(a_{ij}(t)e_j(t)+b_{ij}(t)f_j(t)) \quad \text{or} \quad \Psi=A_tE_t+B_tF_t,
\end{equation}
where $A_t$ is the matrix with entries $a_{ij}(t)$, $B_t$ is the matrix with entries $b_{ij}(t)$, and $\Psi$, $E_t$, and $F_t$ are matrices with rows $dd\mathfrak f(\varsigma_i)$, $e_i(t)$, and $f_i(t)$, respectively.

By a result in \cite{FiRi}, the measure $\varphi_{t*}\eta$ is absolutely continuous with respect to $\eta$. Let $g_t$ be the density of $\varphi_{t*}\eta$ (i.e. $\varphi_{t*}\eta=g_t\eta$). Since $\varphi_t$ is smooth in a neighborhood of $x'$, we can consider $g_t\eta$ as a volume form. If $\{e_i(t),f_j(t)\}$ is a canonical frame at $\alpha$, then $\{e^{s\vec H}_*(e_i(t+s)),e^{s\vec H}_*(f_j(t+s))\}$ is a canonical frame at $e^{s\vec H}(\alpha)$. It follows from this and $\varphi_{t*}\eta=g_t\eta$ that
\begin{equation}\label{gt}
g_t(\varphi_t(x')) \,|\eta(d\varphi_t(\varsigma_1),d\varphi_t(\varsigma_2),d\varphi_t(\varsigma_3))| =|\eta(\varsigma_1,\varsigma_2,\varsigma_3)|.
\end{equation}

On the other hand, it follows from the definition of the canonical Darboux frame, the definition of $\varphi_t$, and (\ref{AB}) that
\[
\begin{split}
d\varphi_t(\varsigma_i)&=d\pi (de^{t\vec H}(dd\mathfrak f(\varsigma_i)))\\
&=\sum\limits_{j=1}^3b_{ij}(t)d\pi (de^{t\vec H}(f_j(t))).
\end{split}
\]
Therefore, this together with Theorem \ref{structuraldetail} gives
\begin{equation}\label{gt2}
|\eta(d\varphi_t(\varsigma_1),d\varphi_t(\varsigma_2),d\varphi_t(\varsigma_3))| =|\eta(\varsigma_1,\varsigma_2,\varsigma_3)\det B_t|.
\end{equation}

Note also that since $(d\varphi_t)_{x'}$ is nonsingular, $B_t$ is invertible for all $t<1$. Since $B_0$ is the identity matrix, $\det B_t>0$ for all $t<1$. Therefore, by combining (\ref{gt}) and (\ref{gt2}), we have the following lemma.

\begin{lem}\label{density}
\[
g_t(\varphi_t(x'))=\frac{1}{\det B_t(x')}.
\]
Recall that $B_t$ is the matrix defined by the canonical frame at $x'$ in (\ref{AB}). Here we write $B_t(x')$ to emphasize its dependence on $x'$.
\end{lem}
If we differentiate (\ref{AB}) with respect to time $t$ and apply Theorem \ref{structural}, then we have
\[
\begin{split}
0&=\dot A_tE_t+A_t\dot E_t+\dot B_tF_t+B_t\dot F_t\\
&=\dot A_tE_t+A_t(C_1E_t+C_2F_t)+\dot B_tF_t-B_t(RE_t+C_1^TF_t),
\end{split}
\]
where
\[
C_1=\left(\begin{array}{ccc}
0 & 0 & 0\\
1 & 0 & 0\\
0 & 0 & 0
\end{array}\right),\quad C_2=\left(\begin{array}{ccc}
1 & 0 & 0\\
0 & 0 & 0\\
0 & 0 & 1
\end{array}\right),
\]
\[
R=\left(\begin{array}{ccc}
r & 0 & 0\\
0 & 0 & 0\\
0 & 0 & 0
\end{array}\right),
\]
\[
r(x')=h_0^2(d\mathfrak f_{x'})+2KH(d\mathfrak f_{x'})=(v_0\mathfrak f(x'))^2-K\mathfrak f(x'),
\]
and $C_1^T$ denotes the transpose of $C_1$.

Therefore, we have the following equations for the matrices $A_t$ and $B_t$.
\begin{equation}\label{AdBd}
\dot A_t+A_tC_1-B_tR=0, \quad \dot B_t+A_tC_2-B_tC_1^T=0.
\end{equation}

If $s_t=g_t(\varphi_t(x))$, then we have, by (\ref{density}) and (\ref{AdBd}), the following:
\[
\det B_t\ddt\det(B_t^{-1})=-tr(B_t^{-1}\dot B_t)=tr(B_t^{-1}A_tC_2).
\]
Therefore, if we let $S_t=B_t^{-1}A_t$, then we have the following lemma.

\begin{lem}\label{densitychange}
\[
\det B_t=e^{-\int_0^t\textbf{tr}(S_sC_2)ds}.
\]
\end{lem}

By (\ref{AdBd}), the matrix $S_t$ defined by $S_t=B_t^{-1}A_t$ satisfies the following matrix Riccati equation
\[
 \dot S_t-R+S_tC_1+C_1^TS_t-S_tC_2S_t=0.
\]

Since $\varphi_1(x)=x_0$ for all $x$, we have $d\varphi_1(\varsigma_i)=0$. Therefore, by Theorem \ref{structuraldetail} and (\ref{AB}), $B_1=0$. Therefore, $S_t^{-1}$ satisfies the following matrix Riccati equation
\begin{equation}\label{Riccati}
\ddt (S_t^{-1})+S_t^{-1}RS_t^{-1}-C_1S_t^{-1}-S_t^{-1}C_1^T+C_2=0\quad \text{and}\quad S^{-1}_1=0.
\end{equation}

Since the coefficient of the above equation does not depend on time $t$, the solution to this equation can be found explicitly by the result in \cite{Le} as follows.

Let us consider the matrix
\[
\mathcal Q=
\left(
\begin{array}{cc}
C_1 & -C_2 \\
R & -C_1^T \\
\end{array}\right)
\]
and the corresponding matrix differential equation $\ddt q=\mathcal Qq$ together with the condition $q(1)=I$.

The fundamental solution is given by
\[
 q(t)=e^{(t-1)\mathcal Q}=\left(
\begin{array}{cccccc}
           \cos\tau_t & 0 & 0 & \frac{\sin\tau_t}{\tau_0}& \frac{1-\cos\tau_t}{\tau_0^2} & 0\\
           -\frac{\sin\tau_t}{\tau_0} & 1 & 0 & \frac{\cos\tau_t-1}{\tau_0^2} & \frac{\sin\tau_t-\tau_t}{\tau_0^3} & 0\\
           0 & 0 & 1 & 0 & 0 & 1-t\\
    -\tau_0\sin\tau_t & 0 & 0 & \cos\tau_t & \frac{\sin\tau_t}{\tau_0} & 0\\
    0 & 0 & 0 & 0 & 1 & 0\\
    0 & 0 & 0 & 0 & 0 & 1
         \end{array}\right).
\]
if $r>0$,

\[
 q(t)=e^{(t-1)\mathcal Q}=\left(
\begin{array}{cccccc}
           1 & 0 & 0 & 1-t& \frac{(1-t)^2}{2} & 0\\
           t-1 & 1 & 0 & -\frac{(1-t)^2}{2} & -\frac{(1-t)^3}{6} & 0\\
           0 & 0 & 1 & 0 & 0 & 1-t\\
    0 & 0 & 0 & 1 & 1-t & 0\\
    0 & 0 & 0 & 0 & 1 & 0\\
    0 & 0 & 0 & 0 & 0 & 1
         \end{array}\right).
\]
if $r=0$,

\[
 q(t)=e^{(t-1)\mathcal Q}=\left(
\begin{array}{cccccc}
           \cosh\tau_t & 0 & 0 & \frac{\sinh\tau_t}{\tau_0}& \frac{\cosh\tau_t-1}{\tau_0^2} & 0\\
           -\frac{\sinh\tau_t}{\tau_0} & 1 & 0 & \frac{1-\cosh\tau_t}{\tau_0^2} & \frac{\tau_t-\sinh\tau_t}{\tau_0^3} & 0\\
           0 & 0 & 1 & 0 & 0 & 1-t\\
    \tau_0\sinh\tau_t & 0 & 0 & \cosh\tau_t & \frac{\sinh\tau_t}{\tau_0} & 0\\
    0 & 0 & 0 & 0 & 1 & 0\\
    0 & 0 & 0 & 0 & 0 & 1
         \end{array}\right).
\]
if $r<0$, where $\tau_t=\sqrt{|r|}(1-t)$.

It follows from \cite[Theorem 1]{Le} that
\[
\begin{array}{ll}
 S_t^{-1}&=\left(
\begin{array}{ccc}
            \frac{\sin\tau_t}{\tau_0}& \frac{1-\cos\tau_t}{\tau_0^2} & 0\\
           \frac{\cos\tau_t-1}{\tau_0^2} & \frac{\sin\tau_t-\tau_t}{\tau_0^3} & 0\\
            0 & 0 & 1-t
     \end{array}\right)\left(
\begin{array}{ccc}
        \cos\tau_t & \frac{\sin\tau_t}{\tau_0} & 0\\
    0 & 1 & 0\\
    0 & 0 & 1
     \end{array}\right)^{-1}\\&
=\left(
\begin{array}{ccc}
           \frac{\tan\tau_t}{\tau_0} & \frac{\cos\tau_t-1}{\tau_0^2\cos\tau_t} & 0 \\
           \frac{\cos\tau_t-1}{\tau_0^2\cos\tau_t} & \frac{\tan\tau_t-\tau_t}{\tau_0^3} & 0 \\
           0 & 0 & 1-t \\
         \end{array}\right),
\end{array}
\]
if $r>0$,
\[
\begin{array}{ll}
 S_t^{-1}&=\left(
\begin{array}{ccc}
            1-t & \frac{(1-t)^2}{2} & 0\\
           -\frac{(1-t)^2}{2} & -\frac{(1-t)^3}{6} & 0\\
            0 & 0 & 1-t
     \end{array}\right)\left(
\begin{array}{ccc}
        1 & 1-t & 0\\
    0 & 1 & 0\\
    0 & 0 & 1
     \end{array}\right)^{-1}\\&
=\left(
\begin{array}{ccc}
           1-t & -\frac{(1-t)^2}{2} & 0 \\
           -\frac{(1-t)^2}{2} & \frac{(1-t)^3}{3} & 0 \\
           0 & 0 & 1-t \\
         \end{array}\right),
\end{array}
\]
if $r=0$, and
\[
 \begin{array}{ll}
 S_t^{-1}&=\left(
\begin{array}{ccc}
            \frac{\sinh\tau_t}{\tau_0}& \frac{\cosh\tau_t-1}{\tau_0^2} & 0\\
           \frac{1-\cosh\tau_t}{\tau_0^2} & \frac{\tau_t-\sinh\tau_t}{\tau_0^3} & 0\\
            0 & 0 & 1-t
     \end{array}\right)\left(
\begin{array}{ccc}
        \cosh\tau_t & \frac{\sinh\tau_t}{\tau_0} & 0\\
    0 & 1 & 0\\
    0 & 0 & 1
     \end{array}\right)^{-1}\\&
=\left(
\begin{array}{ccc}
           \frac{\tanh\tau_t}{\tau_0} & \frac{1-\cosh\tau_t}{\tau_0^2\cosh\tau_t} & 0 \\
           \frac{1-\cosh\tau_t}{\tau_0^2\cosh\tau_t} & \frac{\tau_t-\tanh\tau_t}{\tau_0^3} & 0 \\
           0 & 0 & 1-t \\
         \end{array}\right).
\end{array}
\]
if $r<0$.

Therefore, inverting the above matrix gives the following.
If $r>0$, then
\[
S_t=\left(
\begin{array}{ccc}
           \frac{\tau_0(\sin\tau_t-\tau_t\cos\tau_t)}{\mathcal D} & \frac{\tau_0^2(1-\cos\tau_t)}{\mathcal D} & 0 \\
           \frac{\tau_0^2(1-\cos\tau_t)}{\mathcal D} & \frac{\tau_0^3\sin\tau_t}{\mathcal D} & 0 \\
           0 & 0 & \frac{1}{1-t} \\
         \end{array}\right)
\]
where $\mathcal D=2-2\cos\tau_t-\tau_t\sin\tau_t$.

If $r=0$, then
\[
S_t=\frac{1}{(1-t)^3}\left(
\begin{array}{ccc}
           4(1-t)^2 & 6(1-t) & 0 \\
           6(1-t) & 12 & 0 \\
           0 & 0 & (1-t)^2 \\
         \end{array}\right).
\]

If $r<0$, then
\[
S_t=\left(
\begin{array}{ccc}
\frac{\tau_0(\tau_t\cosh\tau_t-\sinh\tau_t)}{\mathcal D^h}  & \frac{\tau_0^2(\cosh\tau_t-1)}{\mathcal D^h} & 0 \\
           \frac{\tau_0^2(\cosh\tau_t-1)}{\mathcal D^h} & \frac{\tau_0^3\sinh\tau_t}{\mathcal D^h} & 0 \\
           0 & 0 & \frac{1}{1-t} \\
         \end{array}\right)
\]
where $\mathcal D^h=2-2\cosh\tau_t+\tau_t\sinh\tau_t$.

If $r>0$, then
\[
\textbf{tr}(C_2S_t)=\frac{\tau_0(\sin\tau_t-\tau_t\cos\tau_t)}{2-2\cos\tau_t-\tau_t\sin\tau_t}+\frac{1}{1-t}.
\]
If $r=0$, then
\[
\textbf{tr}(C_2S_t)=\frac{5}{1-t}.
\]
If $r<0$, then
\[
\textbf{tr}(C_2S_t)=\frac{\tau_0(\tau_t\cosh\tau_t-\sinh\tau_t)}{2-2\cosh\tau_t+\tau_t\sinh\tau_t}+\frac{1}{1-t}.
\]

If we integrate the above equations, we get
\begin{equation}\label{S11S33g0}
\int_0^t\textbf{tr}(C_2S_s)ds=-\log\left[\frac{(1-t)(2-2\cos\tau_t-\tau_t\sin\tau_t)}{(2-2\cos\tau_0-\tau_0\sin\tau_0)}\right].
\end{equation}
if $r>0$,
\begin{equation}\label{S11S330}
\int_0^t\textbf{tr}(C_2S_s)ds=-\log(1-t)^5.
\end{equation}
if $r=0$, and
\begin{equation}\label{S11S33l0}
\int_0^t\textbf{tr}(C_2S_s)ds=-\log\left[\frac{(1-t)(2-2\cosh\tau_t+\tau_t\sinh\tau_t)}{(2-2\cosh\tau_0+\tau_0\sinh\tau_0)}\right].
\end{equation}
if $r<0$.

Since all the above computations hold for $\eta$-almost all $x'$, we can combine them with Lemma \ref{density} and \ref{densitychange} and obtain
\[
\begin{split}
\eta(\varphi_t(U))&=\int_{\varphi_t(U)}\frac{1}{g_t(x)}d((\varphi_t)_*\eta)(x)\\
&=\int_U\frac{1}{g_t(\varphi_t(x))}d\eta(x)\\
&=\int_U\det B_t\,d\eta(x)\\
&=\int_Ue^{-\int_0^t\textbf{tr}(S_sC_2)ds}d\eta(x)\\
&=\int_U(1-t) \left(\frac{s(\mathfrak k(x),(1-t)D(x))}{s(\mathfrak k(x),D(x))}\right)d\eta(x).
\end{split}
\]
\end{proof}

\smallskip

\section{Proof of Theorem \ref{mainSasa}  and its Consequences}\label{proofmain}

\begin{proof}[Proof of Theorem \ref{mainSasa}]
We use the setup and notations as in the proof of Theorem \ref{SasaVol}. Both Lemma \ref{density} and \ref{densitychange} still hold in this case. The only difference is that the curvature $R_\alpha(t)$ now is not given explicitly and, more importantly, it depends on time $t$.

Let us consider the following matrix Riccati equation with constant coefficients:

\begin{equation}\label{sriccati}
 \ddt ({\tilde S_t}^{-1})+\tilde S_t^{-1}\tilde R\tilde S_t^{-1}-C_1\tilde S_t^{-1}-\tilde S_t^{-1}C_1^T+C_2=0
\end{equation}
together with the condition
\begin{equation}\label{initial}
\tilde S^{-1}_1=0,
\end{equation}
where $\tilde R_\alpha=\left(
\begin{array}{ccc}
           2KH(\alpha) & 0 & 0 \\
           0 & 0 & 0 \\
           0 & 0 & 0 \\
         \end{array}\right)$.

It follows from the assumption of the theorem that
\[
R^{11}_{d\mathfrak f_{x'}}(t)\geq 2KH(d\mathfrak f_{x'})
\]
and
\[
R^{22}_{d\mathfrak f_{x'}}(t)=0.
\]

Note also that solution of (\ref{sriccati}) and (\ref{initial}) are symmetric. Therefore, by comparison theorem of the matrix Riccati equation (see \cite[Theorem 2.1]{FrJaKa}), we have $S_t^{-1}\geq \tilde S_t^{-1}\geq 0$ for $t$ close enough to 1. Here $A\geq B$ means that $A-B$ is nonnegative definite. By monotonicity (see \cite[Proposition V.1.6]{Bh}), $0\leq S_t\leq\tilde S_t$ for $t$ close enough to 1. Since $S_t$ and $\tilde S_t$ also satisfy Riccati equations, we can apply the same comparison principle to $S_t$ and $\tilde S_t$. It follows that $0\leq S_t\leq\tilde S_t$ for all $t$ in $[0,1]$. Therefore,
\[
\textbf{tr}(\tilde S_tC_2)\geq  \textbf{tr} (S_tC_2).
\]

It follows from Lemma \ref{density} and \ref{densitychange} that
\[
g_t(\varphi_t(z))=e^{\int_0^t\textbf{tr} (S_sC_2)ds}\leq e^{\int_0^t\textbf{tr}(\tilde S_sC_2)ds}.
\]

The last term of the above inequality can be computed as in the proof of Theorem \ref{SasaVol} and this finishes one implication.

Conversely, assume that $h_0(\alpha)^2+2H(\alpha)\kappa(x_0)=R^{11}_\alpha=\ric(\alpha)< 2H(\alpha)K$ for some point $x_0$ in the manifold $M$ and some covector $\alpha$ in $T^*_{x_0}M$. By replacing $\alpha$ by $\alpha-h_0(\alpha)\alpha_0$ in the above inequality, we see that $\kappa(x_0)<K$. Let $\delta$ and $\epsilon>0$ be small enough so that $K-\kappa(x)>\epsilon$ for all $x$ inside a subriemannian ball of radius $\delta$.

Let $\mathcal U$ be the set of covectors $\alpha$ in $T^*_{x_0}M$ such that $t\mapsto \pi(e^{t\vec H}(\alpha))$ is minimizing between its end points, $\sqrt{2H(\alpha)}<\delta$, and $h_0(\alpha)^2< 2H(\alpha)\epsilon$. The set $\mathcal U$ has nonzero measure and so is the image $U:=\pi(e^{1\cdot \vec H}(\mathcal U))$. Note that the condition $\sqrt{2H(\alpha)}<\epsilon$ ensures that $\kappa(x)<K-\epsilon$ for all $x$ in $U$. Note also that $-e^{1\vec H}(d\mathfrak f_x)$ is contained in $\mathcal U$. It follows from conservation of $h_0$ and $H$ that
\[
\begin{split}
R^{11}_{d\mathfrak f_x}(t)
&=h_0^2(e^{1\vec H}(d\mathfrak f_x))+2H(e^{1\vec H}(d\mathfrak f_x))\kappa(\varphi_t(x)) \\
&< 2H(e^{1\vec H}(d\mathfrak f_x))(\epsilon+\kappa(\varphi_t(x)))\\
&< 2H(d\mathfrak f_x)K.
\end{split}
\]
for all $x$ in $\mathcal U$.

An argument using comparison theorem of Riccati equation as above shows that the above chosen point $x_0$ and set $U$ violate the definition of $\mathcal{MCP}(K;2,3)$.
\end{proof}

\begin{proof}[Proof of Corollary \ref{doublingSasa}]
From the proof of Theorem \ref{mainSasa}, we have
\[
\eta(\varphi_{1/2}(B_{x_0}(2kR)))\geq  \frac{1}{2^5}\eta(B_{x_0}(2kR)).
\]

Since $\varphi_{1/2}(B_{x_0}(2kR))$ is contained in $B_{x_0}(kR)$, the result follows.
\end{proof}

The proof of Corollary \ref{PoincareSasa} can be found, for instance, in \cite{LoVi2,vo}, for metric spaces satisfying condition $MCP(0,5)$. We give the proof here in the subriemannian case for completeness.

\begin{proof}[Proof of Corollary \ref{PoincareSasa}]
Let $x'$ and $\bar x$ be two points on the manifold $M$. Let $\mathfrak f(x)=-\frac{1}{2}d^2(x,x')$ and let $\bar{\mathfrak f}(x)=-\frac{1}{2}d^2(x,\bar x)$. Let
\[
\varphi_t(x):=\pi(e^{t\vec H}(d\mathfrak f_{x})),\quad \bar\varphi_t(x):=\pi(e^{t\vec H}(d\bar{\mathfrak f}_{x})).
\]

Recall that $t\mapsto\varphi_t(x)$ is a minimizing geodesic connecting $x$ and $x'$ for $\eta$-almost all $x$. Assume that both $x'$ and $\bar x$ are contained in the ball $B_{x_0}(kR)$. Then
\[
|f(\bar x)-f(x')|\leq 2R\left(\int_0^{1/2}|\nabla_H f(\varphi_t(\bar x))|dt +\int_0^{1/2}|\nabla_H f(\bar\varphi_t(x'))|dt\right).
\]

By the proof of Theorem \ref{mainSasa}, we have
\[
\begin{split}
&\int_{B_{x_0}(kR)}\int_0^{1/2}|\nabla_H f(\varphi_t(x))|dt d\eta(x)\\
\leq &\int_0^{1/2}\int_{\varphi_t(B_{x_0}(kR))}\frac{|\nabla_H f(x)|}{(1-t)^5}d\eta(x)dt.
\end{split}
\]

Since $\varphi_t(B_{x_0}(kR))$ is contained in $B_{x_0}(2kR)$, it follows that
\[
\int_{B_{x_0}(kR)}\int_0^{1/2}|\nabla_H f(\varphi_t(x))|dt d\eta(x)
\leq \frac{15}{4}\int_{B_{x_0}(2kR)}|\nabla_H f(x)|d\eta(x).
\]

Therefore,
\[
\begin{split}
&\int_{B_{x_0}(kR)}|f(x')-\left<f\right>_{B_{x_0}(kR)}|d\eta(x')\\
&\leq\frac{1}{\eta(B_{x_0}(kR))}\int_{B_{x_0}(kR)}\int_{B_{x_0}(kR)}|f(\bar x)-f(x')|d\eta(\bar x)d\eta(x')\\
&\leq 15R\int_{B_{x_0}(2kR)}|\nabla_H f(x)|d\eta(x)
\end{split}
\]

Finally, if we apply Corollary \ref{doublingSasa}, then we get
\[
\begin{split}
&\frac{1}{\eta(B_{x_0}(kR))}\int_{B_{x_0}(kR)}|f(x')-\left<f\right>_{B_{x_0}(kR)}|d\eta(x')\\
&\leq \frac{15R}{\eta(B_{x_0}(kR))}\int_{B_{x_0}(2kR)}|\nabla_H f(x)|d\eta(x)\\
&\leq \frac{480R}{\eta(B_{x_0}(2kR))}\int_{B_{x_0}(2kR)}|\nabla_H f(x)|d\eta(x).
\end{split}
\]
\end{proof}

\smallskip

\section{Proof of Theorem \ref{injcompact}}\label{proofinj}

In this section, we give the proof of Theorem \ref{injcompact}. Like the proof of Theorem \ref{R22}, it involves the expansion of $R^{22}$ found after a lengthy calculations. So only a sketch of the proof will be given.

\begin{proof}[Proof of Theorem \ref{injcompact}]
Let $\alpha_i$ be an unbounded sequence in $\Omega$. It follows that $h_0^2(\alpha_i)+H(\alpha_i)\to\infty$ as $i\to\infty$. By \cite[Theorem 3]{Ag2}, it follows that $H(\alpha_i)\to 0$ and so $h_0^2(\alpha_i)\to\infty$ as $i\to\infty$. On the other hand, recall that $R^{22}$ is a degree 4 polynomial of $h_0$, $h_1$, and $h_2$. But  $R^{22}$ is also quadratic in $h_0$ and the coefficient of $R^{22}_{\alpha_i}$ in $h_0^2$ is given by
\[
3(c_{01}^2+c_{02}^1)h_2^2+12c_{01}^1h_1h_2-3(c_{01}^2+c_{02}^1)h_1^2.
\]

Moreover, by \cite[Theorem 3.1]{Ag},
\[
\sqrt{h_1^2(\alpha_i)+h_2^2(\alpha_i)}=\sqrt{2H(\alpha_i)}=\frac{2\pi}{h_0(\alpha_i)}+O\left(\frac{1}{h_0^2(\alpha_i)}\right) \quad \text{as } i\to\infty.
\]

It follows that $R^{22}_{\alpha_i}$ stays bounded as $i\to\infty$ and the bound is independent of the sequence by compactness of the manifold.
\end{proof}

\end{document}